\documentclass{amsart}
\usepackage{amssymb}
\usepackage[all]{xy}

\theoremstyle{remark}
\newtheorem{rem}{Remark}[subsection]
\newtheorem{rems}[rem]{Remarks}
\newtheorem{example}[rem]{Example}
\theoremstyle{definition}
\newtheorem{defn}[rem]{Definition}
\theoremstyle{plain}
\newtheorem{prop}[rem]{Proposition}
\newtheorem{lemma}[rem]{Lemma}
\newtheorem{cor}[rem]{Corollary}
\newtheorem{thm}[rem]{Theorem}

\DeclareMathOperator{\SL}{SL}
\DeclareMathOperator{\GL}{GL}
\DeclareMathOperator{\PGL}{PGL}
\DeclareMathOperator{\Mat}{Mat}

\DeclareMathOperator{\Hom}{Hom}
\DeclareMathOperator{\End}{End}

\DeclareMathOperator{\Mum}{Mum}

\DeclareMathOperator{\id}{id}
\DeclareMathOperator{\pr}{pr}
\DeclareMathOperator{\dt}{dt}

\DeclareMathOperator{\characteristic}{char}
\DeclareMathOperator{\Spec}{Spec}
\DeclareMathOperator{\Pic}{Pic}
\DeclareMathOperator{\NS}{NS}
\DeclareMathOperator{\Lin}{Lin}

\DeclareMathOperator{\Gr}{Gr}
\DeclareMathOperator{\glue}{glue}
\DeclareMathOperator{\diag}{diag}
\DeclareMathOperator{\rank}{rk}
\newcommand{\Bun}{\calM}
\newcommand{\calM}{\mathcal{M}}

\newcommand{\calT}{\mathcal{T}}
\newcommand{\calG}{\mathcal{G}}
\newcommand{\calH}{\mathcal{H}}
\newcommand{\calX}{\mathcal{X}}
\newcommand{\calY}{\mathcal{Y}}
\newcommand{\calL}{\mathcal{L}}

\newcommand{\calE}{\mathcal{E}}
\newcommand{\calPic}{\mathcal{P}\mathit{ic}}
\newcommand{\relPic}{\underline{\Pic}}
\newcommand{\relHom}{\underline{\Hom}}

\newcommand{\M}{\mathfrak{M}}
\newcommand{\longto}[1][]{\stackrel{#1}{\longrightarrow}}

\newcommand{\g}{\mathfrak{g}}

\newcommand{\Gm}{\mathbb{G}_m}
\newcommand{\Q}{\mathbb{Q}}

\newcommand{\Z}{\mathbb{Z}}
\newcommand{\C}{\mathbb{C}}
\newcommand{\coroots}{\mathrm{coroots}}

\newcommand{\univ}{\mathrm{univ}}
\newcommand{\dual}{\mathrm{dual}}
\newcommand{\Luniv}{\mathcal{L}^{\univ}}
\newcommand{\Ldet}{\mathcal{L}^{\det}}
\renewcommand{\O}{\mathcal{O}}
\renewcommand{\H}{\mathrm{H}}

\newcommand{\ab}{\text{ab}}
\newcommand{\ad}{\text{ad}}

\newcommand{\ltilde}{\widetilde{l}}
\newcommand{\Gtilde}{\widetilde{G}}
\newcommand{\Htilde}{\widetilde{H}}
\newcommand{\Lambdatilde}{\widetilde{\Lambda}}
\newcommand{\Ghat}{\widehat{G}}
\newcommand{\Gbar}{\bar{G}}
\newcommand{\Gad}{G^{\ad}}
\newcommand{\Had}{H^{\ad}}

\newcommand{\Hbar}{\bar{H}}

\newcommand{\Ohat}{\widehat{\O}}
\newcommand{\ghat}{\widehat{g}}
\newcommand{\xihat}{\widehat{\xi}}
\newcommand{\pihat}{\widehat{\pi}}

\newcommand{\deltahat}{\widehat{\delta}}
\newcommand{\deltabar}{\bar{\delta}}
\newcommand{\dbar}{\bar{d}}

\begin{document}

\title[Line bundles on moduli stacks of principal bundles]{The line 
bundles on moduli stacks of principal bundles on a curve}

\author[I. Biswas]{Indranil Biswas}
\address{School of Mathematics, Tata Institute of Fundamental Research, Homi Bhabha Road, Mumbai 400005, India}
\email{indranil@math.tifr.res.in}

\author[N. Hoffmann]{Norbert Hoffmann}
\address{Mathematisches Institut der Freien Universit\"at, Arnimallee 3, 14195 Berlin, Germany}
\email{norbert.hoffmann@fu-berlin.de}
\thanks{The second author gratefully acknowledges the support of the SFB/TR 45 "Perioden, Modulr\"aume und Arithmetik algebraischer Variet\"aten"}



\date{}

\begin{abstract}
  Let $G$ be an affine reductive algebraic group over an algebraically 
  closed field $k$. We determine the Picard group of the moduli stacks
  of principal $G$--bundles on any smooth projective curve over $k$.
\end{abstract}

\maketitle

\section{Introduction}
As long as moduli spaces of bundles on a smooth projective algebraic curve
$C$ have been studied, their Picard groups have attracted some interest.
The first case was the coarse moduli scheme of semistable vector bundles 
with fixed determinant over a curve $C$ of genus $g_C \geq 2$. Seshadri
proved that its Picard group is infinite cyclic in the coprime case
\cite{ramanan}; Dr\'{e}zet and Narasimhan showed that this remains
valid in the non--coprime case also \cite{drezet-narasimhan}. 

The case of principal $G$--bundles over $C$ for simply connected, almost 
simple groups $G$ over the complex numbers has been studied intensively,
motivated also by the relation to conformal field theory and the Verlinde
formula \cite{beauville-laszlo, faltings_verlinde, kumar-narasimhan-ramanathan}.
Here Kumar and Narasimhan \cite{kumar-narasimhan} showed that the 
Picard group of the coarse moduli scheme of semistable $G$--principal 
bundles over a curve $C$ of genus $g_C \geq 2$ embeds as a
subgroup of finite index into the Picard group of the affine 
Grassmannian, which is canonically isomorphic to $\mathbb Z$; 
this finite index was determined recently in \cite{boysal-kumar}.
Concerning the Picard group of the moduli stack $\Bun_G$ of principal
$G$--bundles over a curve $C$ of any genus $g_C \geq 0$, Laszlo and 
Sorger \cite{L-S, sorger} showed that its canonical map to the Picard group $\Z$
of the affine Grassmannian is actually an isomorphism. Faltings 
\cite{faltings} has generalised this result to positive characteristic,
and in fact to arbitrary noetherian base scheme.

If $G$ is not simply connected, then the moduli stack $\Bun_G$ has 
several connected components which are indexed by $\pi_1( G)$. For
any $d \,\in\, \pi_1( G)$, let $\Bun_G^d$ be the corresponding
connected component of $\Bun_G$. For semisimple, almost simple
groups $G$ over $\C$, the Picard group $\Pic( \Bun_G^d)$ has been
determined case by case by Beauville, Laszlo and Sorger
\cite{BLS, laszlo}. It is finitely generated, and its torsion part is a 
direct sum of $2 g_C$ copies of $\pi_1( G)$. Furthermore, its 
torsion--free part again embeds as a subgroup of finite index
into the Picard group $\Z$ of the affine Grassmannian. Together
with a general expression for this index, Teleman \cite{teleman}
also proved these statements, using topological and analytic methods.

In this paper, we determine the Picard group $\Pic( \Bun_G^d)$ for any reductive
group $G$, working over an algebraically closed ground field $k$ without
any restriction on the characteristic of $k$ (for all $g_C \geq 0$). Endowing this 
group with a natural scheme structure, we prove that the resulting group scheme
$\relPic( \Bun_G^d)$ over $k$ contains, as an open subgroup, the scheme 
of homomorphisms from $\pi_1( G)$ to the Jacobian $J_C$, with the 
quotient being a finitely generated free abelian group which we
denote by $\NS( \Bun_G^d)$ and call it the N\'{e}ron--Severi group
(see Theorem \ref{mainthm}). We introduce this N\'{e}ron--Severi
group combinatorially in \S~\ref{subsection:NS_reductive};
in particular, Proposition \ref{NS_extension} describes it as
follows: the group $\NS( \Bun_G^d)$ contains a subgroup
$\NS( \Bun_{G^{\ab}})$ which depends only on the torus $G^{\ab} \,=\, 
G/[G\, , G]$; the quotient is a group of Weyl--invariant symmetric 
bilinear forms on the root system of the semisimple part
$[G\, , G]$, subject to certain integrality conditions that
generalise Teleman's result in \cite{teleman}.

We also describe the maps of Picard groups induced by group
homomorphisms $G \longto H$. An interesting effect appears for 
the inclusion $\iota_G: T_G \hookrightarrow G$ of a maximal
torus, say for semisimple $G$: Here the induced map 
$\NS( \Bun_G^d) \longto  \NS( \Bun_{T_G}^{\delta})$ for a lift 
$\delta \in \pi_1( T_G)$ of $d$ involves contracting each
bilinear form in $\NS( \Bun_G^d)$ to a linear form by means
of $\delta$ (cf. Definition \ref{def:iota_G}). In general,
the map of Picard groups induced by a group homomorphism
$G \longto H$ is a combination of this effect and of more
straightforward induced maps (cf. Definition
\ref{NS_reductive:functorial} and Theorem \ref{mainthm}.iv).
In particular, these induced maps are different on different
components of $\Bun_G$, whereas the Picard groups
$\Pic( \Bun_G^d)$ themselves are essentially independent of $d$.

Our proof is based on Faltings' result in the simply connected case. To 
deduce the general case, the strategy of \cite{BLS} and \cite{laszlo} 
is followed, meaning we ``cover'' the moduli stack $\Bun_G^d$ by
a moduli stack of ``twisted'' bundles as in \cite{BLS} under
the universal cover of $G$, more precisely under an
appropriate torus times the universal cover of the semisimple part 
$[G\, , G]$. To this ``covering'', we apply Laszlo's \cite{laszlo}
method of descent for torsors under a group stack. To understand the
relevant descent data, it turns out that we may restrict to the maximal
torus $T_G$ in $G$, roughly speaking because the pullback $\iota_G^*$
is injective on the Picard groups of the moduli stacks.

We briefly describe the structure of this paper. In Section 
\ref{section:Bun_G}, we recall the relevant moduli stacks and collect
some basic facts. Section \ref{section:torus} deals with the case that
$G = T$ is a torus. Section \ref{section:twisted} treats the ``twisted''
simply connected case as indicated above. In the final Section
\ref{section:reductive}, we put everything together to prove our
main theorem, namely Theorem \ref{mainthm}. Each section begins
with a slightly more detailed description of its contents.

Our motivation for this work was to understand the existence of
Poincar\'{e} families on the corresponding coarse moduli schemes,
or in other words to decide whether these moduli stacks are neutral
as gerbes over their coarse moduli schemes. The consequences for this
question will be spelled out in a subsequent paper.

\section{The stack of $G$--bundles and its Picard functor} \label{section:Bun_G}
Here we introduce the basic objects of this paper, namely the moduli 
stack of principal $G$--bundles on an algebraic curve and its Picard functor.
The main purpose of this section is to fix some notation and terminology;
along the way, we record a few basic facts for later use.

\subsection{A Picard functor for algebraic stacks}
Throughout this paper, we work over an algebraically closed field $k$.
There is no restriction on the characteristic of $k$.
We say that a stack $\calX$ over $k$ is \emph{algebraic}
if it is an Artin stack and also locally of finite type over $k$. 
Every algebraic stack $\calX \neq \emptyset$ admits a point
$x_0: \Spec( k) \longto \calX$ according to Hilbert's Nullstellensatz.

A $1$--morphism $\Phi: \calX \longto \calY$ of stacks is an 
\emph{equivalence} if some $1$--morphism $\Psi: \calY \longto \calX$
admits $2$--isomorphisms $\Psi \circ \Phi \cong \id_{\calX}$ and
$\Phi \circ \Psi \cong \id_{\calY}$. A diagram
\begin{equation*} \xymatrix{
  \calX \ar[r]^A \ar[d]_{\Phi} & \calX' \ar[d]^{\Phi'}\\
  \calY \ar[r]^B & \calY'
} \end{equation*}
of stacks and $1$--morphisms is \emph{$2$--commutative} if a 
$2$--isomorphism $\Phi' \circ A \cong B \circ \Phi$ is given. Such a 
$2$--commutative diagram is \emph{$2$--cartesian} if the induced
$1$--morphism from $\calX$ to the fibre product of stacks
$\calX' \times_{\calY'} \calY$ is an equivalence.

Let $\calX$ and $\calY$ be algebraic stacks over $k$. As usual, we
denote by $\Pic( \calX)$ the abelian group of isomorphism classes of
line bundles $\calL$ on $\calX$. If $\calX \neq \emptyset$, then
\begin{equation*}
  \pr_2^*: \Pic( \calY) \longto \Pic( \calX \times \calY)
\end{equation*}
is injective because $x_0^*: \Pic( \calX \times \calY) \longto \Pic( \calY)$
is a left inverse of $\pr_2^*$.
\begin{defn} \label{def:Pic}
  The \emph{Picard functor} $\relPic( \calX)$ is the functor from
  schemes $S$ of finite type over $k$ to abelian groups that
  sends $S$ to $\Pic( \calX \times S)/\pr_2^* \Pic( S)$.
\end{defn}
If $\relPic( \calX)$ is representable, then we denote the representing scheme
again by $\relPic( \calX)$. If $\relPic( \calX)$ is the constant sheaf given by
an abelian group $\Lambda$, then we say that $\relPic( \calX)$ is
\emph{discrete} and simply write $\relPic( \calX) \cong \Lambda$.
(Since the constant Zariski sheaf $\Lambda$ is already an fppf sheaf,
it is not necessary to specify the topology here.)
\begin{lemma} \label{lemma:Pic}
  Let $\calX$ and $\calY$ be algebraic stacks over $k$ with $\Gamma( \calX, \O_{\calX}) = k$.
  \begin{itemize}
   \item[i)] The canonical map
    \begin{equation*}
      \pr_2^*: \Gamma( \calY, \O_{\calY}) \longto \Gamma( \calX \times \calY, \O_{\calX \times \calY})
    \end{equation*}
    is an isomorphism.
   \item[ii)] Let $\calL \in \Pic( \calX \times \calY)$ be given.
    If there is an atlas $u: U \longto \calY$ for which
    $u^* \calL \in \Pic( \calX \times U)$ is trivial,
    then $\calL \in \pr_2^* \Pic( \calY)$.
  \end{itemize}
\end{lemma}
\begin{proof}
  i) Since the question is local in $\calY$, we may assume that
  $\calY = \Spec( A)$ is an affine scheme over $k$. In this case, we have
  \begin{equation*}
    \Gamma( \calX \times \calY, \O_{\calX \times \calY}) = \Gamma( \calX, (\pr_1)_* \O_{\calX \times \calY})
      = \Gamma( \calX, A \otimes_k \O_{\calX}) = A = \Gamma( \calY, \O_{\calY}).
  \end{equation*}

  ii) Choose a point $x_0: \Spec( k) \longto \calX$.
  We claim that $\calL$ is isomorphic to $\pr_2^* \calL_{x_0}$
  for $\calL_{x_0} := x_0^* \calL \in \Pic( \calY)$. More precisely
  there is a unique isomorphism $\calL \cong \pr_2^* \calL_{x_0}$
  whose restriction to $\{x_0\} \times \calY \cong \calY$ is the
  identity. To prove this, due to the uniqueness involved, this
  claim is local in $\calY$. Hence we may assume $\calY = U$, which
  by assumption means that $\calL$ is trivial.
  In this case, statement (i) implies the claim.
\end{proof}
\begin{cor} \label{Pic:injective}
  For $\nu = 1, 2$, let $\calX_{\nu}$ be an algebraic stack
  over $k$ with $\Gamma( \calX_{\nu}, \O_{\calX_{\nu}}) = k$.
  Let $\Phi: \calX_1 \longto \calX_2$ be a $1$--morphism
  such that the induced morphism of functors
  $\Phi^*: \relPic( \calX_2) \longto \relPic( \calX_1)$ is injective. Then
  \begin{equation*}
    \Phi^*: \Pic( \calX_2 \times \calY) \longto \Pic( \calX_1 \times \calY)
  \end{equation*}
  is injective for every algebraic stack $\calY$ over $k$.
\end{cor}
\begin{proof}
  Since $\calY$ is assumed to be locally of finite type over $k$, we can 
  choose an atlas $u: U \longto \calY$ such that $U$ is a disjoint union 
  of schemes of finite type over $k$. Suppose that
  $\calL \in \Pic( \calX_2 \times \calY)$ has trivial pullback
  $\Phi^* \calL \in \Pic( \calX_1 \times \calY)$.
  Then $(\Phi \times u)^* \calL \in \Pic( \calX_1 \times U)$ is also 
  trivial. Using the assumption on $\Phi^*$ it follows that 
  $u^* \calL \in \Pic( \calX_2 \times U)$ is trivial.
  Now apply Lemma \ref{lemma:Pic}(ii).
\end{proof}

We will also need the following stacky version of the standard see--saw 
principle.
\begin{lemma} \label{see-saw}
  Let $\calX$ and $\calY$ be two nonempty algebraic stacks over $k$. If
  $\relPic( \calX)$ is discrete, and $\Gamma( \calY, \O_{\calY}) = k$, then
  \begin{equation*}
    \pr_1^* \oplus \pr_2^*: \relPic( \calX) \oplus \relPic( \calY) \longto \relPic( \calX \times \calY)
  \end{equation*}
  is an isomorphism of functors.
\end{lemma}
\begin{proof}
  Choose points $x_0: \Spec( k) \longto \calX$ and
  $y_0: \Spec( k) \longto \calY$. The morphism of functors
  $\pr_1^* \oplus \pr_2^*$ in question is injective, because
  \begin{equation*}
    y_0^* \oplus x_0^*: \relPic( \calX \times \calY) \longto \relPic( \calX) \oplus \relPic( \calY)
  \end{equation*}
  is a left inverse of it. Therefore, to prove the lemma it suffices to 
  show that $y_0^* \oplus x_0^*$ is also injective.

  So let a scheme $S$ of finite type over $k$ be given,
  as well as a line bundle $\calL$ on $\calX \times \calY \times S$ such 
  that $y_0^* \calL$ is trivial in $\relPic( \calX)$. We claim that $\calL$
  is isomorphic to the pullback of a line bundle on $\calY \times S$.
 
  To prove the claim, tensoring $\calL$ with an appropriate line bundle
  on $S$ if necessary, we may assume that $y_0^* \calL$ is trivial in
  $\Pic( \calX \times S)$. By assumption, $\relPic( \calX) \cong \Lambda$
  for some abelian group $\Lambda$. Sending any
  $(y, s): \Spec( k) \longto \calY \times S$ to the isomorphism class of
  \begin{equation*}
    (y, s)^*( \calL) \in \Pic( \calX)
  \end{equation*}
  we obtain a Zariski--locally constant map from the set of $k$--points
  in $\calY \times S$ to $\Lambda$. This map vanishes on $\{y_0\} \times S$,
  and hence it vanishes identically on $\calY \times S$ because $\calY$ is 
  connected. This means that $u^* \calL \in \Pic( \calX \times U)$
  is trivial for any atlas $u: U \longto \calY \times S$.
  Now Lemma \ref{lemma:Pic}(ii) completes the proof of the claim.

  If moreover $x_0^* \calL$ is trivial in $\relPic( \calY)$, then $\calL$ is even isomorphic to the
  pullback of a line bundle on $S$, and hence trivial
  in $\relPic(\calX \times \calY)$. This proves the injectivity of 
  $y_0^* \oplus x_0^*$, and hence the lemma follows.
\end{proof}

\subsection{Principal $G$--bundles over a curve}
We fix an irreducible smooth projective curve $C$ over the 
algebraically closed base field $k$. The genus of $C$ will be
denoted by $g_C$. Given a linear algebraic group
$G \hookrightarrow \GL_n$, we denote by
\begin{equation*}
  \Bun_G
\end{equation*}
the moduli stack of principal $G$--bundles $E$ on $C$. More precisely, 
$\Bun_G$ is given by the groupoid $\Bun_G( S)$ of principal $G$--bundles 
on $S \times C$ for every $k$--scheme $S$. The stack $\Bun_G$ is known
to be algebraic over $k$ (see for example \cite[Proposition 3.4]{L-S}, or
\cite[Th\'{e}or\`{e}me 4.6.2.1]{L-MB} together with
\cite[Lemma 4.8.1]{ramanathan_II}).

Given a morphism of linear algebraic groups $\varphi: G \longto H$, 
the extension of the structure group by $\varphi$ defines
a canonical $1$--morphism
\begin{equation*}
  \varphi_*: \Bun_G \longto \Bun_H
\end{equation*}
which more precisely sends a principal $G$--bundle $E$
to the principal $H$--bundle
\begin{equation*}
  \varphi_* E := E \times^G H := (E \times G)/H,
\end{equation*}
following the convention that principal bundles carry a
\emph{right} group action. One has a canonical $2$--isomorphism
$(\psi \circ \varphi)_* \cong \psi_* \circ \varphi_*$ whenever
$\psi: H \longto K$ is another morphism of linear algebraic groups.
\begin{lemma} \label{stacks:cartesian}
  Suppose that the diagram of linear algebraic groups
  \begin{equation*} \xymatrix{
    H \ar[r]^{\psi_2} \ar[d]_{\psi_1} & G_2 \ar[d]^{\varphi_2}\\
    G_1 \ar[r]_{\varphi_1} & G
  } \end{equation*}
  is cartesian. Then the induced $2$--commutative diagram
  of moduli stacks
  \begin{equation*} \xymatrix{
    \Bun_H \ar[r]^{(\psi_2)_*} \ar[d]_{(\psi_1)_*} & \Bun_{G_2} \ar[d]^{(\varphi_2)_*}\\
    \Bun_{G_1} \ar[r]_{(\varphi_1)_*} & \Bun_G
  } \end{equation*}
  is $2$--cartesian.
\end{lemma}
\begin{proof}
  The above $2$--commutative diagram defines a $1$--morphism
  \begin{equation*}
    \Bun_H \longto \Bun_{G_1} \times_{\Bun_G} \Bun_{G_2}.
  \end{equation*}
  To construct an inverse, let $E$ be a principal $G$--bundle on some 
  $k$--scheme $X$. For $\nu = 1, 2$, let $E_{\nu}$ be a principal
  $G_{\nu}$--bundle on $X$ together with an isomorphism
  $E_{\nu} \times^{G_{\nu}} G \cong E$; note that the latter defines
  a map $E_{\nu} \longto E$ of schemes over $X$. Then $G_1 \times G_2$
  acts on $E_1 \times_X E_2$, and the closed subgroup
  $H \subseteq G_1 \times G_2$ preserves the closed subscheme
  \begin{equation*}
    F := E_1 \times_E E_2 \subseteq E_1 \times_X E_2\, .
  \end{equation*}
  This action turns $F$ into a principal $H$--bundle.
  Thus we obtain in particular a $1$--morphism
  \begin{equation*}
    \Bun_{G_1} \times_{\Bun_G} \Bun_{G_2} \longto \Bun_H.
  \end{equation*}
  It is easy to check that this is the required inverse.
\end{proof}
Let $Z$ be a closed subgroup in the center of $G$. Then the 
multiplication $Z \times G \longto G$ is a group
homomorphism; we denote the induced $1$--morphism by
\begin{equation*}
  \_ \otimes \_: \Bun_Z \times \Bun_G \longto \Bun_G
\end{equation*}
and call it \textit{tensor product}. In particular,
tensoring with a principal $Z$--bundle $\xi$ on $C$
defines a $1$--morphism which we denote by
\begin{equation} \label{txi}
  t_{\xi}: \Bun_G \longto \Bun_G.
\end{equation}
For commutative $G$, this tensor product makes $\Bun_G$ a group stack.

Suppose now that $G$ is reductive. We follow the convention that all reductive groups are smooth and connected.
In particular, $\Bun_G$ is also smooth \cite[4.5.1]{behrend}, so its connected components and its
irreducible components coincide; we denote this set
of components by $\pi_0( \Bun_G)$. This set
$\pi_0( \Bun_G)$ can be described as follows:

Let $\iota_G: T_G \hookrightarrow G$ be the inclusion of a maximal torus, with cocharacter group
$\Lambda_{T_G} := \Hom( \Gm, T_G)$. Let $\Lambda_{\coroots} \subseteq \Lambda_{T_G}$ be the subgroup
generated by the coroots of $G$. The Weyl group of $(G, T_G)$ acts 
trivially on $\Lambda_{T_G}/\Lambda_{\coroots}$, so this quotient is,
up to a \emph{canonical} isomorphism, independent of the choice of $T_G$.
We denote this quotient by $\pi_1( G)$; if $\pi_1( G)$ is trivial,
then $G$ is called simply connected. For $k = \C$, these definitions
coincide with the usual notions for the topological space $G( \C)$. 

Sending each line bundle on $C$ to its degree we define an isomorphism 
$\pi_0( \Bun_{\Gm}) \longto \Z$, which induces an isomorphism
$\pi_0( \Bun_{T_G}) \longto \Lambda_{T_G}$. Its inverse, composed with the map
\begin{equation*}
  (\iota_G)_*: \pi_0( \Bun_{T_G}) \longto \pi_0( \Bun_G)\, ,
\end{equation*}
is known to induce a canonical bijection
\begin{equation*}
  \pi_1( G) = \Lambda_{T_G}/\Lambda_{\coroots} \longto[ \sim] \pi_0( \Bun_G), 
\end{equation*}
cf. \cite{D-S} and \cite{holla}. We denote by $\Bun^d_G$ the component of $\Bun_G$ given by $d \in \pi_1( G)$.
\begin{lemma} \label{flat}
  Let $\varphi: G \twoheadrightarrow H$ be an epimorphism of reductive groups over $k$ whose kernel is contained in the center of $G$.
  For each $d \in \pi_1( G)$, the $1$--morphism
  \begin{equation*}
    \varphi_*: \Bun_G^d \longto \Bun_H^e, \qquad e := \varphi_*( d) \in \pi_1( H),
  \end{equation*}
  is faithfully flat.
\end{lemma}
\begin{proof}
  Let $T_H \subseteq H$ be the image of the maximal torus
  $T_G \subseteq G$. Let $B_G \subseteq G$ be a Borel
  subgroup containing $T_G$; then
  \begin{equation*}
    B_H\, :=\, \varphi(B_G)\, \subset\, H
  \end{equation*}
  is a Borel subgroup of $H$ containing $T_H$.
  For the moment, we denote
  \begin{itemize}
   \item by $\Bun_{T_G}^d \subseteq \Bun_{T_G}$ and $\Bun_{B_G}^d \subseteq \Bun_{B_G}$ the inverse images of $\Bun^d_G \subseteq \Bun_G$, and
   \item by $\Bun_{T_H}^e \subseteq \Bun_{T_H}$ and $\Bun_{B_H}^e \subseteq \Bun_{B_H}$ the inverse images of $\Bun^e_H \subseteq \Bun_H$.
  \end{itemize}
  Let $\pi_G: B_G \twoheadrightarrow T_G$ and $\pi_H: B_H \twoheadrightarrow T_H$ denote the canonical surjections. Then
  \begin{equation*}
    \Bun_{B_G}^d = (\pi_G)_*^{-1}( \Bun_{T_G}^d) \qquad\text{and}\qquad \Bun_{B_H}^e = (\pi_H)_*^{-1}( \Bun_{T_H}^e),
  \end{equation*}
  because $\pi_0( \Bun_{T_G}) = \pi_0( \Bun_{B_G})$ and
  $\pi_0( \Bun_{T_H}) = \pi_0( \Bun_{B_H})$ according to
  the proof of \cite[Proposition 5]{D-S}. Applying Lemma 
  \ref{stacks:cartesian} to the two cartesian squares
  \begin{equation*} \xymatrix{
    T_G \ar[d]_{\varphi_T} & B_G \ar[d]^{\varphi_B} \ar@{^{(}->}[r] \ar@{->>}[l]_{\pi_G} & G \ar[d]^{\varphi}\\
    T_H & B_H \ar@{^{(}->}[r] \ar@{->>}[l]_{\pi_H} & H
  } \end{equation*}
  of groups, we get two $2$--cartesian squares
  \begin{equation*} \xymatrix{
    \Bun_{T_G}^d \ar[d]_{(\varphi_T)_*} & \Bun_{B_G}^d \ar[l] \ar[r] \ar[d]^{(\varphi_B)_*} & \Bun_G^d \ar[d]^{\varphi_*}\\
    \Bun_{T_H}^e & \Bun_{B_H}^e \ar[l] \ar[r] & \Bun_H^e\\
  } \end{equation*}
  of moduli stacks. Since $(\varphi_T)_*$ is faithfully flat, its pullback $(\varphi_B)_*$ is so as well.
  This implies that $\varphi_*$ is also faithfully flat, as some open substack of $\Bun_{B_H}^e$ maps
  smoothly and surjectively onto $\Bun_H^e$, according to \cite[Propositions 1 and 2]{D-S}.
\end{proof}

\section{The case of torus} \label{section:torus}
This section deals with the Picard functor of the moduli stack 
$\Bun_G^0$ in the special case where $G = T$ is a torus. We explain in 
the second subsection that its description involves the character group
$\Hom( T, \Gm)$ and the Picard functor of its coarse moduli scheme,
which is isomorphic to a product of copies of the Jacobian $J_C$.
As a preparation, the first subsection deals with the N\'{e}ron--Severi 
group of such products of principally polarised abelian varieties.
A little care is required to keep everything functorial in $T$, since 
this functoriality will be needed later.

\subsection{On principally polarised abelian varieties}

Let $A$ be an abelian variety over $k$, with dual abelian variety 
$A^{\vee}$ and N\'{e}ron--Severi group
\begin{equation*}
  \NS( A) := \Pic( A)/A^{\vee}( k)\, .
\end{equation*}
For a line bundle $L$ on $A$, the standard morphism
\begin{equation*}
  \phi_L: A \longto A^{\vee}
\end{equation*}
sends $a \in A( k)$ to $\tau_a^*( L) \otimes L^{\dual}$ where
$\tau_a: A \longto A$ is the translation by $a$.
$\phi_L$ is a homomorphism by the theorem of the cube
\cite[\S 6]{mumford}. Let a principal polarisation
\begin{equation*}
  \phi: A \longto[ \sim] A^{\vee}
\end{equation*}
be given. Let
\begin{equation*}
  c^{\phi}: \NS( A) \longto \End A
\end{equation*}
be the injective homomorphism that sends the class
of $L$ to $\phi^{-1} \circ \phi_L$. We denote by
$\dagger: \End A \longto \End A$ the Rosati involution associated
to $\phi$; so by definition, it sends $\alpha: A \longto A$
to $\alpha^{\dagger} := \phi^{-1} \circ \alpha^{\vee} \circ \phi$.

\begin{lemma} \label{lemma:ppav}
  An endomorphism $\alpha \in \End( A)$ is in the image of $c^{\phi}$
  if and only if $\alpha^{\dagger} = \alpha$.
\end{lemma}
\begin{proof}
  If $k = \C$, this is contained in
  \cite[Chapter 5, Proposition 2.1]{lange-birkenhake}.
  For polarisations of arbitrary degree, the analogous statement
  about $\End( A) \otimes \Q$ is shown in \cite[p. 190]{mumford};
  its proof carries over to the situation of this lemma as follows.

  Let $l$ be a prime number, $l \neq \characteristic( k)$, and let
  \begin{equation*}
    e_l: T_l( A) \times T_l( A^{\vee}) \longto \Z_l( 1)
  \end{equation*}
  be the standard pairing between the Tate modules of $A$ and 
  $A^{\vee}$, cf. \cite[\S 20]{mumford}. According to
  \cite[\S 20, Theorem 2 and \S 23, Theorem 3]{mumford}, a homomorphism
  $\psi: A \longto A^{\vee}$ is of the form $\psi = \phi_L$
  for some line bundle $L$ on $A$ if and only if
  \begin{equation*}
    e_l( x, \psi_* y) = - e_l( y, \psi_* x)
      \quad\text{for all}\quad x, y \in T_l( A)\, .
  \end{equation*}
  In particular, this holds for $\phi$. Hence the right hand side equals
  \begin{equation*}
    - e_l( y, \psi_* x) = - e_l( y, \phi_* \phi^{-1}_* \psi_* x) = e_l( \phi^{-1}_* \psi_* x, \phi_* y)
    = e_l( x, \psi^{\vee}_* (\phi^{-1})^{\vee}_* \phi_* y)\, ,
  \end{equation*}
  where the last equality follows from \cite[p. 186, equation (I)]{mumford}. Since the pairing $e_l$ is nondegenerate, it
  follows that $\psi = \phi_L$ holds for some $L$ if and only if
  \begin{equation*}
    \psi_* y = \psi^{\vee}_* (\phi^{-1})^{\vee}_* \phi_* y 
      \quad\text{for all}\quad y \in T_l( A)\, ,
  \end{equation*}
  hence if and only if $\psi = \psi^{\vee} \circ (\phi^{-1})^{\vee} \circ \phi$.
  By definition of the Rosati involution $\dagger$, the latter is
  equivalent to $(\phi^{-1} \circ \psi)^{\dagger} = \phi^{-1} \circ \psi$.
\end{proof}
Let $\Lambda$ be a finitely generated free abelian group.
Let $\Lambda \otimes A$ denote the abelian variety over $k$
with group of $S$--valued points $\Lambda \otimes A( S)$
for any $k$--scheme $S$.
\begin{defn}
  The subgroup
  \begin{equation*}
    \Hom^s( \Lambda \otimes \Lambda, \End A) \subseteq \Hom( \Lambda \otimes \Lambda, \End A)
  \end{equation*}
  consists of all $b: \Lambda \otimes \Lambda \longto \End A$ with
  $b( \lambda_1 \otimes \lambda_2)^{\dagger} = b( \lambda_2 \otimes \lambda_1)$
  for all $\lambda_1, \lambda_2 \in \Lambda$.
\end{defn}
\begin{cor} \label{cor:ppav}
  There is a unique isomorphism
  \begin{equation*}
    c_{\Lambda}^{\phi}: \NS( \Lambda \otimes A) \longto[ \sim] \Hom^s( \Lambda \otimes \Lambda, \End A)
  \end{equation*}
  which sends the class of each line bundle $L$ on $\Lambda \otimes A$ to the linear map
  \begin{equation*}
    c_{\Lambda}^{\phi}( L): \Lambda \otimes \Lambda \longto \End A
  \end{equation*}
  defined by sending $\lambda_1 \otimes \lambda_2$ for $\lambda_1, \lambda_2 \in \Lambda$ to the composition
  \begin{equation*}
    A \xrightarrow{ \lambda_1 \otimes \_} \Lambda \otimes A \longto[ \phi_L] (\Lambda \otimes A)^{\vee}
      \xrightarrow{ (\lambda_2 \otimes \_)^{\vee}} A^{\vee}  \longto[ \phi^{-1}] A.
  \end{equation*}
\end{cor}
\begin{proof}
  The uniqueness is clear. For the existence, we may then choose an isomorphism
  $\Lambda \cong \Z^r$; it yields an isomorphism $\Lambda \otimes A \cong A^r$. Let
  \begin{equation*}
    \phi^r = \underbrace{\phi \times \cdots \times \phi}_{r \text{ factors}}: A^r \longto[ \sim] (A^{\vee})^r = (A^r)^{\vee}
  \end{equation*}
  be the diagonal principal polarisation on $A^r$.
  According to Lemma \ref{lemma:ppav},
  \begin{equation*}
    c^{\phi^r}: \NS( A^r) \longto \End( A^r)
  \end{equation*}
  is an isomorphism onto the Rosati--invariants.
  Under the standard isomorphisms
  \begin{equation*}
    \End( A^r) = \Mat_{r \times r}( \End A) = \Hom( \Z^r \otimes \Z^r, \End A),
  \end{equation*}
  the Rosati involution on $\End( A^r)$ corresponds to the involution
  $( \alpha_{ij}) \longmapsto ( \alpha^{\dagger}_{ji})$ on
  $\Mat_{r \times r}( \End A)$, and hence the Rosati--invariant part
  of $\End( A^r)$ corresponds to $\Hom^s( \Z^r \otimes \Z^r, \End A)$.
  Thus we obtain an isomorphism
  \begin{equation*}
    \NS( \Lambda \otimes A) \cong \NS( A^r) \longto[ c^{\phi^r}] \Hom^s( \Z^r \otimes \Z^r, \End A) \cong \Hom^s( \Lambda \otimes \Lambda, \End A).
  \end{equation*}
  By construction, it maps the class of each line
  bundle $L$ on $\Lambda \otimes A$ to the map
  $c_{\Lambda}^{\phi}( L): \Lambda \otimes \Lambda \longto \End A$
  prescribed above.
\end{proof}

\subsection{Line bundles on $\Bun_T^0$} \label{subsection:torus}
Let $T \cong \Gm^r$ be a torus over $k$. We will always denote by
\begin{equation*}
  \Lambda_T := \Hom( \Gm, T)
\end{equation*}
the cocharacter lattice. We set in the previous subsection this 
finitely generated free abelian group and the Jacobian variety $J_C$,
endowed with the standard principal polarisation
$\phi_{\Theta}: J_C \longto[ \sim] J_C^{\vee}$.
\begin{defn} \label{NS:torus}
  The finitely generated free abelian group
  \begin{equation*}
    \NS( \Bun_T) := \Hom( \Lambda_T, \Z) \oplus \Hom^s( \Lambda_T \otimes \Lambda_T, \End J_C)
  \end{equation*}
  is the \emph{N\'{e}ron--Severi group} of $\Bun_T$.
\end{defn}
For each finitely generated abelian group $\Lambda$, we denote by 
$\relHom( \Lambda, J_C)$ the $k$--scheme of homomorphisms from $\Lambda$ to $J_C$.
If $\Lambda \cong \Z^r \times \Z/n_1 \times \cdots \times \Z/n_s$, then
\begin{equation*}
  \relHom( \Lambda, J_C) \cong J_C^r \times J_C[ n_1] \times \cdots \times J_C[ n_s]
\end{equation*}
where $J_C[ n]$ denotes the kernel of the map
$J_C \longto J_C$ defined by multiplication with $n$.
\begin{prop} \label{prop:torus} \begin{itemize}
 \item[i)] The Picard functor $\relPic( \Bun_T^0)$ is representable by a scheme locally of finite type over $k$.
 \item[ii)] There is a canonical exact sequence of commutative group schemes
  \begin{equation*}
    0 \longto \relHom( \Lambda_T, J_C) \longto[j_T] \relPic( \Bun_T^0) \longto[ c_T] \NS( \Bun_T) \longto 0.
  \end{equation*}
 \item[iii)] Let $\xi$ be a principal $T$--bundle of
  degree $0 \in \Lambda_T$ on $C$. Then the diagram
  \begin{equation*} \xymatrix{
    0 \ar[r] & \relHom( \Lambda_T, J_C) \ar[r]^-{j_T} \ar@{=}[d] & \relPic( \Bun_T^0) \ar[r]^{c_T} \ar[d]^{t_{\xi}^*} & \NS( \Bun_T) \ar[r] \ar@{=}[d] & 0\\
    0 \ar[r] & \relHom( \Lambda_T, J_C) \ar[r]^-{j_T} & \relPic( \Bun_T^0) \ar[r]^{c_T} & \NS( \Bun_T) \ar[r] & 0
  } \end{equation*}
  commutes.
\end{itemize} \end{prop}
\begin{proof}
  Given a line bundle $\calL$ on $\Bun^0_T$, the automorphism
  group $T$ of each point in $\Bun^0_T$ acts on the corresponding
  fibre of $\calL$, so we obtain a character
  \begin{equation*}
    w( \calL): T \longto \Gm
  \end{equation*}
  which is independent of the choice of the point in $\Bun^0_T$. This
  character is called the \emph{weight} $w( \calL)$ of $\calL$. Let
  \begin{equation*}
    q: \Bun_T^0 \longto \M_T^0
  \end{equation*}
  be the canonical morphism to the coarse moduli scheme $\M_T^0$,
  which is an abelian variety canonically isomorphic to
  $\relHom( \Lambda_T, J_C)$. Line bundles of weight $0$
  on $\Bun_T^0$ descend to $\M_T^0$, so the sequence
  \begin{equation*}
    0 \longto \Pic( \M_T^0) \longto[ q^*] \Pic( \Bun_T^0) \longto[ w] \Hom( \Lambda_T, \Z)
  \end{equation*}
  is exact. This extends for families. Since $\relPic( A)$
  is representable for any abelian variety $A$,
  the proof of (i) is now complete.

  Standard theory of abelian varieties and Corollary
  \ref{cor:ppav} together yield another short exact sequence
  \begin{equation*}
    0 \longto \relHom( \Lambda_T, J_C) \longto \relPic( \M_T^0)
      \longto \Hom^s( \Lambda_T \otimes \Lambda_T, \End J_C) \longto 0\, .
  \end{equation*}
  Given a character $\chi: T \longto \Gm$ and $p \in C( k)$, we denote 
  by $\chi_* \Luniv_p$ the line bundle on $\Bun_T^0$ that
  associates to each $T$--bundle $L$ on $C$ the $\Gm$--bundle
  $\chi_* L_p$. Clearly, $\chi_* \Luniv_p$ has weight $\chi$; in
  particular, it follows that $w$ is surjective, so we get
  an exact sequence of discrete abelian groups
  \begin{equation*}
    0 \longto \Hom^s( \Lambda_T \otimes \Lambda_T, \End J_C)
      \longto \relPic( \Bun_T^0)/\relHom( \Lambda_T, J_C)
      \longto \Hom( \Lambda_T, \Z) \longto 0\, .
  \end{equation*}
  Since $C$ is connected, the algebraic equivalence class of
  $\chi_* \Luniv_p$ does not depend on the choice of $p$;
  sending $\chi$ to the class of $\chi_* \Luniv_p$ thus defines a
  canonical splitting of the latter exact sequence. This proves (ii).

  Finally, it is standard that $t_{\xi}^*$ (see \eqref{txi}) is the 
  identity map on $\relPic^0( \M_T^0) = \relHom( \Lambda, J_C)$
  (see \cite[Proposition 9.2]{milne}), and $t_{\xi}^*$ induces the
  identity map on the discrete quotient $\relPic( \Bun_T^0)/\relPic^0( \M_T^0)$
  because $\xi$ can be connected to the trivial $T$--bundle in $\Bun_T^0$.
\end{proof}

\begin{rem} \label{c_T:functorial} \upshape
  The exact sequence in Proposition \ref{prop:torus}(ii) is functorial 
  in $T$. More precisely, each homomorphism of tori $\varphi: T \longto T'$
  induces a morphism of exact sequences
  \begin{equation*} \xymatrix{
    0 \ar[r] & \relHom( \Lambda_{T'}, J_C) \ar[r]^-{j_{T'}} \ar[d]^{\varphi^*}
             & \relPic( \Bun_{T'}^0) \ar[r]^{c_{T'}} \ar[d]^{\varphi^*} & \NS( \Bun_{T'}) \ar[r] \ar[d]^{\varphi^*} & 0\\
    0 \ar[r] & \relHom( \Lambda_T, J_C) \ar[r]^-{j_T} & \relPic( \Bun_T^0) \ar[r]^{c_T} & \NS( \Bun_T) \ar[r] & 0.
  } \end{equation*}
\end{rem}
\begin{cor} \label{torus:product}
  Let $T_1$ and $T_2$ be tori over $k$. Then
  \begin{equation*}
    \pr_1^* \oplus \pr_2^*: \relPic( \Bun^0_{T_1}) \oplus \relPic( \Bun^0_{T_2}) \longto \relPic( \Bun^0_{T_1 \times T_2})
  \end{equation*}
  is a closed immersion of commutative group schemes over $k$.
\end{cor}
\begin{proof}
  As before, let $\Lambda_{T_1}$, $\Lambda_{T_2}$ and
  $\Lambda_{T_1 \times T_2}$ denote the cocharacter lattices. Then
  \begin{equation*}
    \pr_1^* \oplus \pr_2^*: \relHom( \Lambda_{T_1}, J_C) \oplus \relHom( \Lambda_{T_2}, J_C) \longto \relHom( \Lambda_{T_1 \times T_2}, J_C)
  \end{equation*}
  is an isomorphism, and the homomorphism of discrete abelian groups
  \begin{equation*}
    \pr_1^* \oplus \pr_2^*: \NS( \Bun_{T_1}) \oplus \NS( \Bun_{T_2}) \longto \NS( \Bun_{T_1 \times T_2})
  \end{equation*}
  is injective by Definition \ref{NS:torus}.
\end{proof}

\section{The twisted simply connected case} \label{section:twisted}
Throughout most of this section, the reductive group $G$ over $k$ will be simply
connected. Using the work of Faltings \cite{faltings} on the Picard functor
of $\Bun_G$, we describe here the Picard functor of the twisted moduli 
stacks $\Bun_{\Ghat, L}$ introduced in \cite{BLS}. In the case $G = \SL_n$, these are
moduli stacks of vector bundles with fixed determinant; their construction
in general is recalled in Subsection \ref{subsection:Grass} below.

The result, proved in that subsection as Proposition \ref{prop:grass}, 
is essentially the same: for almost simple $G$, line bundles on $\Bun_{\Ghat, L}$
are classified by an integer, their so--called central charge. The main 
tool for that are as usual algebraic loop groups; what we need about them
is collected in Subsection \ref{loopgroups}.

For later use, we need to keep track of the functoriality in $G$, in particular
of the pullback to a maximal torus $T_G$ in $G$. To state this more easily, we
translate the central charge into a Weyl--invariant symmetric bilinear 
form on the cocharacter lattice of $T_G$, replacing each integer by the
corresponding multiple of the basic inner product. This allows to describe the 
pullback to $T_G$ in Proposition \ref{prop:twisted}(iii). Along the way, 
we also consider the pullback along representations of $G$; these
just correspond to the pullback of bilinear forms, which reformulates
--- and generalises to arbitrary characteristic --- the usual
multiplication by the Dynkin index \cite{kumar-narasimhan-ramanathan}.
Subsection \ref{sc:NS} describes these pullback maps combinatorially
in terms of the root system, and Subsection \ref{sc2torus} proves that
these combinatorial maps actually give the pullback of line bundles
on these moduli stacks.

\subsection{Loop groups} \label{loopgroups}
Let $G$ be a reductive group over $k$. We denote
\begin{itemize}
 \item by $LG$ the algebraic loop group of $G$, meaning the group 
  ind--scheme over $k$ whose group of $A$--valued points for any 
  $k$--algebra $A$ is $G( A(( t)))$,
 \item by $L^+ G \subseteq LG$ the subgroup with $A$--valued points
  $G( A[[ t]]) \subseteq G( A((t)))$,
 \item and for $n \geq 1$, by $L^{\geq n} G \subseteq L^+ G$
  the kernel of the reduction modulo $t^n$.
\end{itemize}
Note that $L^+ G$ and $L^{\geq n} G$ are affine group schemes over $k$.
The $k$-algebra corresponding to $L^{\geq n} G$ is the inductive limit
over all $N > n$ of the $k$--algebras corresponding to $L^{\geq n} G/L^{\geq N}$.
A similar statement holds for $L^+ G$.

If $X$ is anything defined over $k$, let $X_S$ denote
its pullback to a $k$--scheme $S$.
\begin{lemma}
  Let $S$ be a reduced scheme over $k$. For $n \geq 1$,
  every morphism $\varphi: (L^{\geq n} G)_S \longto (\Gm)_S$
  of group schemes over $S$ is trivial.
\end{lemma}
\begin{proof}
  This follows from the fact that $L^{\geq n} G$
  is pro--unipotent; more precisely:

  As $S$ is reduced, the claim can be checked on geometric points $\Spec( k') \longto S$.
  Replacing $k$ by the larger algebraically closed field $k'$ if necessary, we may thus
  assume $S = \Spec( k)$; then $\varphi$ is a morphism $L^{\geq n} G \longto \Gm$.
  
  Since the $k$--algebra corresponding to $\Gm$ is finitely generated, it follows
  that $\varphi$ factors through $L^{\geq n} G/L^{\geq N}$
  for some $N > n$. Denoting by $\g$ the Lie algebra of $G$,
  \cite[II, \S4, Theorem 3.5]{D-G} provides an exact sequence
  \begin{equation*}
    1 \longto L^{\geq N} G \longto L^{\geq N-1} G \longto \g \longto 1.
  \end{equation*}
  Here $\varphi$ restricts to a character of $\g$, which has to vanish; thus $\varphi$ also
  factors through $L^{\geq n} G/L^{\geq N-1}$. Iterating this argument shows that $\varphi$ is trivial.
\end{proof}
\begin{lemma} \label{split}
  Suppose that the reductive group $G$ is simply connected, in particular semisimple.
  If a central extension of group schemes over $k$
  \begin{equation} \label{central_ext}
    1 \longto \Gm \longto \calH \longto[ \pi] L^+ G \longto 1
  \end{equation}
  splits over $L^{\geq n} G$ for some $n \geq 1$, then it splits over $L^+ G$.
\end{lemma}
\begin{proof}
  Let a splitting over $L^{\geq n} G$ be given, i.\,e.\ a homomorphism of 
  group schemes $\sigma: L^{\geq n} G \longto \calH$ such that
  $\pi \circ \sigma = \id$. Given points $h \in \calH( S)$ and
  $g \in L^{\geq n} G( S)$ for some $k$-scheme $S$, the two elements
  \begin{equation*}
    h \cdot \sigma( g) \cdot h^{-1} \qquad\text{and}\qquad \sigma( \pi( h) \cdot g \cdot \pi( h)^{-1})
  \end{equation*}
  in $\calH( S)$ have the same image under $\pi$, so their difference is an element in $\Gm( S)$,
  which we denote by $\varphi_h( g)$. Sending $h$ and $g$ to $h$ and $\varphi_h( g)$ defines
  a morphism
  \begin{equation*}
    \varphi: (L^{\geq n} G)_{\calH} \longto (\Gm)_{\calH}
  \end{equation*}
  of group schemes over $\calH$. Since $L^+ G/L^{\geq 1} G \cong G$ and $L^{\geq N-1} G/L^{\geq N} \cong \g$
  for $N \geq 2$ are smooth, their successive extension $L^+ G/L^{\geq N} G$ is also smooth.
  Thus the limit $L^+ G$ is reduced, so $\calH$ is reduced as well.
  Using the previous lemma, it follows that $\varphi$ is the constant map $1$; in other
  words, $\sigma$ commutes with conjugation. $\sigma$ is a closed immersion because $\pi \circ \sigma$ is, so
  $\sigma$ is an isomorphism onto a closed normal subgroup, and the quotient is a central extension
  \begin{equation*}
    1 \longto \Gm \longto \calH \big/ \sigma( L^{\geq n} G) \longto L^+ G \big/ L^{\geq n} G \longto 1.
  \end{equation*}
  If $n \geq 2$, then this restricts to a central extension of
  $L^{\geq n - 1} G/L^{\geq n} G \cong \g$ by $\Gm$. It can be
  shown that any such extension splits.

  (Indeed, the unipotent radical of the extension projects
  isomorphically to the quotient $\g$. Note that the unipotent
  radical does not intersect the subgroup $\Gm$, and 
  the quotient by the subgroup generated by the unipotent radical
  and $\Gm$ is reductive, so this this reductive quotient being a
  quotient of $\g$ is in fact trivial.)

  Therefore, the image of a section $\g \longto \calH \big/ \sigma( L^{\geq n} G)$
  has an inverse image in $\calH$ which $\pi$ maps isomorphically onto $L^{\geq n - 1} G \subseteq L^+ G$.
  Hence the given central extension \eqref{central_ext} splits over $L^{\geq n - 1} G$ as well. Repeating this argument,
  we get a splitting over $L^{\geq 1} G$, and finally also over $L^+ G$,
  because every central extension of $L^+ G/L^{\geq 1} G \cong G$
  by $\Gm$ splits as well, $G$ being simply connected.

  (To prove the last assertion, for any extension $\Gtilde$
  of $G$ by $\Gm$, consider the commutator subgroup
  $[\Gtilde\, ,\Gtilde]$ of $\Gtilde$.
  It projects surjectively to the commutator subgroup of
  $G$ which is $G$ itself. Since $[\Gtilde\, ,\Gtilde]$
  is connected and reduced, and $G$ is simply connected,
  this surjective morphism must be an isomorphism.)
\end{proof}

\subsection{Descent from the affine Grassmannian} \label{subsection:Grass}
Let $G$ be a reductive group over $k$. We denote by $\Gr_G$ the affine 
Grassmannian of $G$, i.\,e.\ the quotient $LG/L^+ G$ in the category of fppf--sheaves. 
Given a point $p \in C( k)$ and a uniformising element $z \in \Ohat_{C, p}$,
there is a standard $1$--morphism
\begin{equation*}
  \glue_{p, z}: \Gr_G \longto \Bun_G
\end{equation*}
that sends each coset $f \cdot L^+ G$ to the trivial $G$--bundles
over $C \setminus \{p\}$ and over $\Ohat_{C, p}$, glued by the
automorphism $f( z)$ of the trivial $G$--bundle over the
intersection; cf. for example \cite[Section 3]{L-S},
\cite[Corollary 16]{faltings}, or \cite[Proposition 3]{jochen}.

For the rest of this subsection, we assume that $G$ is simply connected, 
hence semisimple. In this case, $\Gr_G$ is known to be an ind--scheme over $k$.
More precisely, \cite[Theorem 8]{faltings} implies that $\Gr_G$ is an inductive
limit of projective Schubert varieties over $k$, which are reduced and
irreducible. Thus the canonical map
\begin{equation} \label{fcts_Gr}
  \pr_2^*: \Gamma( S, \O_S) \longto \Gamma( \Gr_G \times S, \O_{\Gr_G \times S})
\end{equation}
is an isomorphism for every scheme $S$ of finite type over $k$.

Define the Picard functor $\relPic( \Gr_G)$ from schemes of finite type
over $k$ to abelian groups as in definition \ref{def:Pic}. The following
theorem about it is proved in full generality in \cite{faltings}. Over
$k = \C$, the group $\Pic( \Gr_G)$ is also determined in \cite{mathieu}
as well as in \cite{kumar-narasimhan-ramanathan}, and $\Pic( \Bun_G)$
is determined in  \cite{L-S} together with \cite{sorger}.
\begin{thm}[Faltings] \label{faltings}
  Let $G$ be simply connected and almost simple.
  \begin{itemize}
   \item[i)] $\relPic( \Gr_G) \cong \Z$.
   \item[ii)] $\glue_{p, z}^*: \relPic( \Bun_G) \longto \relPic( \Gr_G)$ 
    is an isomorphism of functors.
  \end{itemize}
\end{thm}
The purpose of this subsection is to carry part (ii) over to twisted 
moduli stacks in the sense of \cite{BLS}; cf. also the first remark
on page 67 of \cite{faltings}. More precisely, let an exact sequence
of reductive groups
\begin{equation} \label{dtr}
  1 \longto G \longto \Ghat \longto[ \dt] \Gm \longto 1
\end{equation}
be given, and a line bundle $L$ on $C$. We denote by $\Bun_{\Ghat, L}$ 
the moduli stack of principal $\Ghat$--bundles $E$ on $C$ together
with an isomorphism $\dt_* E \cong L$; cf. section 2 of \cite{BLS}. If 
for example the given exact sequence is
\begin{equation*}
  1 \longto \SL_n \longto \GL_n \longto[ \det] \Gm \longto 1,
\end{equation*}
then $\Bun_{\GL_n, L}$ is the moduli stack of vector bundles with fixed determinant $L$.

In general, the stack $\Bun_{\Ghat, L}$ comes with a $2$--cartesian diagram
\begin{equation*} \xymatrix{
  \Bun_{\Ghat, L} \ar[r] \ar[d] & \Bun_{\Ghat} \ar[d]^{\dt_*}\\\Spec( k) \ar[r]^L & \Bun_{\Gm}
} \end{equation*}
from which we see in particular that $\Bun_{\Ghat, L}$ is algebraic.
It satisfies the following variant of the Drinfeld--Simpson 
uniformisation theorem \cite[Theorem 3]{D-S}.
\begin{lemma} \label{uniformise}
  Let a point $p \in C( k)$ and a principal $\Ghat$--bundle $\calE$ on 
  $C \times S$ for some $k$--scheme $S$ be given. Every trivialisation
  of the line bundle $\dt_* \calE$ over $(C \setminus \{p\}) \times S$
  can \'{e}tale--locally in $S$ be lifted to a trivialisation
  of $\calE$ over $(C \setminus \{p\}) \times S$.
\end{lemma}
\begin{proof}
  The proof in \cite{D-S} carries over to this situation as follows.
  Choose a maximal torus $T_{\Ghat} \subseteq \Ghat$. Using
  \cite[Theorem 1]{D-S}, we may assume that $\calE$ comes from a 
  principal $T_{\Ghat}$--bundle; cf. the first paragraph in the
  proof of \cite[Theorem 3]{D-S}. Arguing as in the third paragraph of 
  that proof, we may change this principal $T_{\Ghat}$--bundle by the
  extension of $\Gm$--bundles along coroots $\Gm \longto T_{\Ghat}$. 
  Since simple coroots freely generate the kernel $T_G$ of
  $T_{\Ghat} \twoheadrightarrow \Gm$, we can thus achieve that this
  $T_{\Ghat}$--bundle is trivial over $(C \setminus \{p\}) \times S$.
  Because $\Gm$ is a direct factor of $T_{\Ghat}$, we can hence lift the given
  trivialisation to the $T_{\Ghat}$--bundle, and hence also to $\calE$.
\end{proof}

Let $d \in \Z$ be the degree of $L$. Since $\dt$ in \eqref{dtr} maps
the (reduced) identity component $Z^0 \cong \Gm$ of the center in
$\Ghat$ surjectively onto $\Gm$, there is a $Z^0$--bundle $\xi$
(of degree $0$) on $C$ with $\dt_*( \xi) \otimes \O_C( dp) \cong L$;
tensoring with it defines an equivalence
\begin{equation*}
  t_{\xi}: \Bun_{\Ghat, \O_C( dp)} \longto[ \sim] \Bun_{\Ghat, L}.
\end{equation*}
Choose a homomorphism $\delta: \Gm \longto \Ghat$ with
$\dt \circ \delta = d \in \Z = \Hom( \Gm, \Gm)$. We denote by
$t^{\delta} \in L\Ghat( k)$ the image of the tautological loop
$t \in L\Gm( k)$ under $\delta_*: L\Gm \longto L\Ghat$. The map
\begin{equation*}
  t^{\delta} \cdot \_: \Gr_G \longto \Gr_{\Ghat}
\end{equation*}
sends, for each point $f$ in $LG$, the coset $f \cdot L^+ G$
to the coset $t^{\delta} f \cdot L^+ \Ghat$. Its composition
$\Gr_G \longto \Bun_{\Ghat}$ with $\glue_{p, z}$ factors
naturally through a $1$--morphism
\begin{equation*}
  \glue_{p, z, \delta}: \Gr_G \longto \Bun_{\Ghat, \O_C( dp)},
\end{equation*}
because $\dt_* \circ ( t^{\delta} \cdot \_): LG \longto L\Ghat \longto L\Gm$ 
is the constant map $t^d$, which via gluing yields the line bundle $\O_C( dp)$.
Lemma \ref{uniformise} provides local sections of $\glue_{p, z, \delta}$.
These show in particular that
\begin{equation*}
  \glue_{p, z, \delta}^*: \Gamma( \Bun_{\Ghat, \O_C( dp)}, \O_{\Bun_{\Ghat, \O_C( dp)}}) \longto \Gamma( \Gr_G, \O_{\Gr_G})
\end{equation*}
is injective. Hence both spaces of sections contain only the constants,
since $\Gamma( \Gr_G, \O_{\Gr_G}) = k$ by equation \eqref{fcts_Gr}.
Using the above equivalence $t_{\xi}$, this implies
\begin{equation} \label{fcts_simple}
  \Gamma( \Bun_{\Ghat, L}, \O_{\Bun_{\Ghat, L}}) = k\, .
\end{equation}
\begin{prop} \label{prop:grass}
  Let $G$ be simply connected and almost simple. Then
  \begin{equation*}
    \glue_{p, z, \delta}^*: \relPic( \Bun_{\Ghat, \O_C( dp)}) \longto \relPic( \Gr_G)
  \end{equation*}
  is an isomorphism of functors.
\end{prop}
\begin{proof}
  $LG$ acts on $\Gr_G$ by multiplication from the left. Embedding the 
  $k$--algebra $\O_{C \setminus p} := \Gamma(C \setminus \{p\}, \O_C)$ into
  $k(( t))$ via the Laurent development at $p$ in the variable $t = z$,
  we denote by $L_{C \setminus p} G \subseteq L G$ the subgroup with
  $A$--valued points $G( A \otimes_k \O_{C \setminus p}) \subseteq G( A(( t)))$
  for any $k$--algebra $A$. The map $\glue_{p, z}$ is a torsor under
  $L_{C \setminus p} G$; cf. for example \cite[Theorem 1.3]{L-S} or
  \cite[Corollary 16]{faltings}. More generally, Lemma \ref{uniformise}
  implies that $\glue_{p, z, \delta}$ is a torsor under the conjugate
  \begin{equation*}
    L_{C \setminus p}^{\delta} G:= t^{-\delta} \cdot L_{C \setminus p} G \cdot t^{\delta} \subseteq L \Ghat
  \end{equation*}
  (which is actually contained in $LG$ since $LG$ is normal in $L \Ghat$), 
  because the action of $L_{C \setminus p}^{\delta} G$ corresponds to
  changing trivialisations over $C \setminus \{p\}$.

  Let $S$ be a scheme of finite type over $k$. Each line bundle on
  $S \times \Bun_{\Ghat, L}$ with trivial pullback to $S \times \Gr_G$
  comes from a character $(L_{C \setminus p}^{\delta} G)_S \longto (\Gm)_S$,
  since the map \eqref{fcts_Gr} is bijective. But
  $L_{C \setminus p}^{\delta} G$ is isomorphic to $L_{C \setminus p} G$,
  and every character $(L_{C \setminus p} G)_S \longto (\Gm)_S$ is 
  trivial according to \cite[p. 66f.]{faltings}. This already shows that
  the morphism of Picard functors $\glue_{p, z, \delta}^*$ is injective.

  The action of $LG$ on $\Gr_G$ induces the trivial action on
  $\relPic( \Gr_G) \cong \Z$, for example because it preserves ampleness,
  or alternatively because $LG$ is connected. Let a line bundle $\calL$ 
  on $\Gr_G$ be given. We denote by $\Mum_{LG}( \calL)$ the Mumford group.
  So $\Mum_{LG}( \calL)$ is the functor from schemes of finite type over 
  $k$ to groups that sends $S$ to the group of pairs $(f, g)$
  consisting of an element $f \in LG( S)$ and an isomorphism
  $g: f^* \calL_S \longto[ \sim] \calL_S$ of line bundles on $\Gr_G \times S$.

  If $f = 1$, then $g \in \Gm( S)$ due to the bijectivity of \eqref{fcts_Gr},
  while for arbitrary $f \in LG( S)$, the line bundles $\calL_S$ and
  $f^* \calL_S$ have the same image in $\relPic( \Gr_G)( S)$,
  implying that $\calL_S$ and $f^* \calL_S$ are Zariski--locally in $S$
  isomorphic. Consequently, we have a short exact sequence
  of sheaves in the Zariski topology
  \begin{equation} \label{extension0}
    1 \longto \Gm \longto \Mum_{LG}( \calL) \longto[ q] LG \longto 1\, .
  \end{equation}
  This central extension splits over $L^+ G \subseteq LG$, because
  the restricted action of $L^+ G$ on $\Gr_G$ has a fixed point. We have to
  show that it also splits over $L_{C \setminus p}^{\delta} G \subseteq LG$.

  Note that $L_{C \setminus p}^{\delta} G = \gamma( L_{C \setminus p} G)$ for
  the automorphism $\gamma$ of $LG$ given by conjugation with $t^{\delta}$.
  Hence it is equivalent to show that the central extension
  \begin{equation} \label{extension}
    1 \longto \Gm \longto \Mum_{LG}( \calL) \xrightarrow{ \gamma^{-1} \circ q} LG \longto 1
  \end{equation}
  splits over $L_{C \setminus p} G$. We know already that it splits over
  $\gamma^{-1}( L^+ G)$, in particular over $L^{\geq n} G$ for some $n \geq 1$.
  Thus it also splits over $L^+ G$, due to Lemma \ref{split}. Hence it 
  comes from a line bundle on $LG/L^+ G = \Gr_G$ (whose associated $\Gm$--bundle
  has total space $\Mum_{LG}( \calL)/L^+ G$, where $L^+ G$ acts from the 
  right via the splitting). According to Theorem \ref{faltings}(ii),
  this line bundle admits a $L_{C \setminus p} G$--linearisation, and 
  hence the extension \eqref{extension} splits indeed over $L_{C \setminus p} G$.

  Thus the extension \eqref{extension0} splits over $L_{C \setminus p}^{\delta} G$,
  so $\calL$ admits an $L_{C \setminus p}^{\delta} G$--linearisation and
  consequently descends to $\Bun_{\Ghat, \O_C( dp)}$. This proves that
  $\glue_{p, z, \delta}^*$ is surjective as a homomorphism of Picard groups.
  Hence it is also surjective as a morphism of Picard functors, because 
  $\relPic( \Gr_G) \cong \Z$ is discrete by Theorem \ref{faltings}(i).
\end{proof}

\begin{rem} \label{glue_functorial}
  Put $\Gad := G/Z$, where $Z \subseteq G$ denotes the center. Given a 
  representation $\rho: \Gad \longto \SL( V)$, we denote its compositions 
  with the canonical epimorphisms $G \twoheadrightarrow \Gad$ and
  $\Ghat \twoheadrightarrow \Gad$ also by $\rho$. Then the diagram
  \begin{equation*} \xymatrix{
    \relPic( \Bun_{\SL( V)}) \ar[d]_{\rho^*} \ar[rr]^-{\glue_{p, z}^*}     && \relPic( \Gr_{\SL( V)}) \ar[d]^{\rho^*}\\
    \relPic( \Bun_{\Ghat, L})        \ar[rr]^-{(t_{\xi} \circ \glue_{p, z, \delta})^*} && \relPic( \Gr_G)
  } \end{equation*}
  commutes.
\end{rem}
\begin{proof}
  Let $t^{\rho \circ \delta} \in L \SL( V)$ denote the image
  of the canonical loop $t \in L \Gm$ under the composition
  $\rho \circ \delta: \Gm \longto \SL( V)$. Then the left
  part of the diagram
  \begin{equation*} \xymatrix{
    & \Gr_G \ar[rr]^-{\glue_{p, z, \delta}} \ar[dl]_{\rho_*} \ar[d]^{t^{\delta} \cdot \_}
      && \Bun_{\Ghat, \O_C( dp)} \ar[d] \ar[r]^{t_{\xi}} & \Bun_{\Ghat, L} \ar[d]\\
    \Gr_{\SL( V)} \ar[rd]_{t^{\rho \circ \delta} \cdot \_} & \Gr_{\Ghat} \ar[rr]^-{\glue_{p, z}} \ar[d]^{\rho_*}
      && \Bun_{\Ghat} \ar[d]_{\rho_*} \ar[r]^{t_{\xi}} & \Bun_{\Ghat} \ar[d]^{\rho_*}\\
     & \Gr_{\SL( V)} \ar[rr]^-{\glue_{p, z}} && \Bun_{\SL( V)} \ar@{=}[r] & \Bun_{\SL( V)}
  } \end{equation*}
  commutes. The four remaining squares are $2$--commutative by 
  construction of the $1$--morphisms $\glue_{p, z, \delta}$,
  $\glue_{p, z}$ and $t_{\xi}$. Applying $\relPic$ to the exterior
  pentagon yields the required commutative square, as $L \SL( V)$
  acts trivially on $\relPic( \Gr_{\SL( V)})$.
\end{proof}

\subsection{N\'{e}ron--Severi groups $\NS( \Bun_G)$ for simply connected $G$} \label{sc:NS}
Let $G$ be a reductive group over $k$; later in this subsection,
we will assume that $G$ is simply connected. Choose a maximal
torus $T_G \subseteq G$, and let
\begin{equation} \label{Hom^W}
  \Hom( \Lambda_{T_G} \otimes \Lambda_{T_G}, \Z)^W
\end{equation}
denote the abelian group of bilinear forms
$b: \Lambda_{T_G} \otimes \Lambda_{T_G} \longto \Z$
that are invariant under the Weyl group
$W = W_G$ of $(G, T_G)$.

Up to a \emph{canonical} isomorphism, the group \eqref{Hom^W} does not 
depend on the choice of $T_G$. More precisely, let $T_{G}' \subseteq G$ be
another maximal torus; then the conjugation $\gamma_g: G \longto G$ with
some $g \in G( k)$ provides an isomorphism from $T_G$ to $T_G'$, and the
induced isomorphism from $\Hom( \Lambda_{T_G'} \otimes \Lambda_{T_G'}, \Z)^W$
to $\Hom( \Lambda_{T_G} \otimes \Lambda_{T_G}, \Z)^W$ does not depend on the 
choice of $g$.

The group \eqref{Hom^W} is also functorial in $G$. More precisely, let 
$\varphi: G \longto H$ be a homomorphism of reductive groups over $k$.
Choose a maximal torus $T_H \subseteq H$ containing $\varphi( T_G)$.
\begin{lemma} \label{from_Weyl}
  Let $T_G' \subseteq G$ be another maximal torus, and let $T_H' \subseteq H$ be a maximal torus containing $\varphi( T_G')$.
  For every $g \in G( k)$ with $T_G' = \gamma_g( T_G)$, there is an $h \in H( k)$ with $T_H' = \gamma_h( T_H)$ such that the following diagram commutes:
  \begin{equation*} \xymatrix{
    T_G \ar[r]^{\gamma_g} \ar[d]^{\varphi} & T_G' \ar[d]^{\varphi}\\ T_H \ar[r]^{\gamma_h} & T_H'
  } \end{equation*}
\end{lemma}
\begin{proof}
  The diagram
  \begin{equation*} \xymatrix{
    T_G \ar@{=}[r] \ar[d]^{\varphi} & T_G \ar[r]^{\gamma_g} \ar[d]^{\varphi} & T_G' \ar[d]^{\varphi}\\
    T_H \ar@{.>}[r] & \gamma_{\varphi( g)}^{-1}( T_H') \ar[r]^-{\gamma_{\varphi( g)}} & T_H'
  } \end{equation*}
  allows us to assume $T_G' = T_G$ and $g = 1$ without loss of generality. Then
  $T_H$ and $T_H'$ are maximal tori in the centraliser of $\varphi( T_G)$,
  which is reductive according to \cite[26.2. Corollary A]{humphreys}. 
  Thus $T_H' = \gamma_h( T_H)$ for an appropriate $k$--point $h$ of this centraliser,
  and $\gamma_h \circ \varphi = \varphi$ on $T_G$ by definition of the centraliser.
\end{proof}
Applying the lemma with $T_G' = T_G$ and $T_H' = T_H$, we see that the 
pullback along $\varphi_*: \Lambda_{T_G} \longto \Lambda_{T_H}$
of a $W_H$--invariant form $\Lambda_{T_H} \otimes \Lambda_{T_H} \longto \Z$ 
is $W_G$--invariant, so we get an induced map
\begin{equation} \label{map:Hom^W}
  \varphi^*: \Hom( \Lambda_{T_H} \otimes \Lambda_{T_H}, \Z)^{W_H} \longto \Hom( \Lambda_{T_G} \otimes \Lambda_{T_G}, \Z)^{W_G}
\end{equation}
which does not depend on the choice of $T_G$ and $T_H$ by the above lemma again.

\emph{For the rest of this subsection, we assume that $G$ and $H$ are simply connected.}
\begin{defn} \label{NS:simply_connected}
  \begin{itemize}
   \item[i)] The \emph{N\'{e}ron--Severi group} $\NS( \Bun_G)$ is the subgroup
    \begin{equation*}
      \NS( \Bun_G) \subseteq \Hom( \Lambda_{T_G} \otimes \Lambda_{T_G}, \Z)^W
    \end{equation*}
    of symmetric forms $b: \Lambda_{T_G} \otimes \Lambda_{T_G} \longto \Z$ 
    with $b( \lambda \otimes \lambda) \in 2\Z$ for all $\lambda \in \Lambda_{T_G}$.
   \item[ii)] Given a homomorphism $\varphi: G \longto H$, we denote by
    \begin{equation*}
      \varphi^*: \NS( \Bun_H) \longto \NS( \Bun_G)
    \end{equation*}
    the restriction of the induced map $\varphi^*$ in \eqref{map:Hom^W}.
  \end{itemize}
\end{defn}
\begin{rems} \label{simply_connected_products}
  i) If $G = G_1 \times G_2$ for simply connected groups $G_1$ and $G_2$, then
  \begin{equation*}
    \NS( \Bun_G) = \NS( \Bun_{G_1}) \oplus \NS( \Bun_{G_2}),
  \end{equation*}
  since each element of $\Hom( \Lambda_{T_G} \otimes \Lambda_{T_G}, \Z)^{W_G}$ vanishes
  on $\Lambda_{T_{G_1}} \otimes \Lambda_{T_{G_2}} + \Lambda_{T_{G_2}} \otimes \Lambda_{T_{G_1}}$.

  ii) If on the other hand $G$ is almost simple, then
  \begin{equation*}
    \NS( \Bun_G) = \Z \cdot b_G
  \end{equation*}
  where the \emph{basic inner product} $b_G$ is the unique element of
  $\NS( \Bun_G)$ that satisfies $b_G( \alpha^{\vee}, \alpha^{\vee}) = 2$
  for all short coroots $\alpha^{\vee} \in \Lambda_{T_G}$ of $G$.

  iii) Let $G$ and $H$ be almost simple. The \emph{Dynkin index} 
  $d_{\varphi} \in \Z$ of a homomorphism $\varphi: G \longto H$
  is defined by $\varphi^*( b_H) = d_{\varphi} \cdot b_G$, cf. \cite[\S 2]{dynkin}.
  If $\varphi$ is nontrivial, then $d_{\varphi} > 0$, since $b_G$ and $b_H$ are positive definite.
\end{rems}
Let $Z \subseteq G$ be the center. Then $\Gad := G/Z$ contains $T_{\Gad} := T_G/Z$ as a maximal torus,
with cocharacter lattice $\Lambda_{T_{\Gad}} \subseteq \Lambda_{T_G} \otimes \Q$.

We say that a homomorphism $l: \Lambda \longto \Lambda'$ between finitely 
generated free abelian groups $\Lambda$ and $\Lambda'$ is \emph{integral}
on a subgroup $\Lambdatilde \subseteq \Lambda \otimes \Q$ if its
restriction to $\Lambda \cap \Lambdatilde$ admits a linear extension
$\ltilde: \Lambdatilde \longto \Lambda'$. By abuse of language,
we will not distinguish between $l$ and its unique linear extension $\ltilde$.
\begin{lemma} \label{extend_b}
  Every element $b: \Lambda_{T_G} \otimes \Lambda_{T_G} \longto \Z$ of 
  $\NS( \Bun_G)$ is integral on $\Lambda_{T_{\Gad}} \otimes \Lambda_{T_G}$
  and on $\Lambda_{T_G} \otimes \Lambda_{T_{\Gad}}$.
\end{lemma}
\begin{proof}
  Let $\alpha: \Lambda_{T_G} \otimes \Q \longto \Q$ be a root of $G$, 
  with corresponding coroot $\alpha^{\vee} \in \Lambda_{T_G}$.
  Lemme 2 in \cite[Chapitre VI, \S 1]{bourbaki_lie} implies the formula
  \begin{equation*}
    b( \lambda \otimes \alpha^{\vee}) = \alpha( \lambda) \cdot b( \alpha^{\vee} \otimes \alpha^{\vee})/2
  \end{equation*}
  for all $\lambda \in \Lambda_{T_G}$. Thus
  $b( \_ \otimes \alpha^{\vee}): \Lambda_{T_G} \longto \Z$
  is an integer multiple of $\alpha$; hence it is integral on
  $\Lambda_{T_{\Gad}}$, the largest subgroup of $\Lambda_{T_G} \otimes \Q$
  on which all roots are integral. But the coroots $\alpha^{\vee}$ generate
  $\Lambda_{T_G}$, as $G$ is simply connected.
\end{proof}
Now let $\iota_G: T_G \hookrightarrow G$ denote the inclusion of the chosen maximal torus.
\begin{defn} \label{def:iota_G}
  Given $\delta \in \Lambda_{T_{\Gad}}$, the homomorphism
  \begin{equation*}
    (\iota_G)^{\NS, \delta}: \NS( \Bun_G) \longto \NS( \Bun_{T_G})
  \end{equation*}
  sends $b: \Lambda_{T_G} \otimes \Lambda_{T_G} \longto \Z$ to
  \begin{equation*}
    b( -\delta \otimes \_): \Lambda_{T_G} \longto \Z \qquad\text{and}\qquad
      \id_{J_C} \cdot b: \Lambda_{T_G} \otimes \Lambda_{T_G} \longto \End J_C.
  \end{equation*}
\end{defn}
This map $(\iota_G)^{\NS, \delta}$ is injective if $g_C \geq 1$,
because all multiples of $\id_{J_C}$ are then nonzero in $\End J_C$.
If $g_C = 0$, then $\End J_C = 0$, but we still have the following
\begin{lemma} \label{injective}
  Every coset $d \in \Lambda_{T_{\Gad}}/\Lambda_{T_G} = \pi_1( \Gad)$
  admits a lift $\delta \in \Lambda_{T_{\Gad}}$ such that the map
  $(\iota_G)^{\NS, \delta}: \NS( \Bun_G) \longto \NS( \Bun_{T_G})$ is injective.
\end{lemma}
\begin{proof}
  Using Remark \ref{simply_connected_products}, we may assume that $G$ 
  is almost simple. In this case, $(\iota_G)^{\NS, \delta}$ is injective
  whenever $\delta \neq 0$, because $\NS( \Bun_G)$ is cyclic and its
  generator $b_G: \Lambda_{T_G} \otimes \Lambda_{T_G} \longto \Z$ is
  as a bilinear form nondegenerate.
\end{proof}
\begin{rem} \label{iotaG_functorial} \upshape
  Given $\varphi: G \longto H$, let $\iota_H: T_H \hookrightarrow H$ be a maximal
  torus with $\varphi( T_G) \subseteq T_H$. If $\delta \in \Lambda_{T_G}$, or if
  more generally $\delta \in \Lambda_{T_{\Gad}}$ is mapped to $\Lambda_{\Had}$ by
  $\varphi_*: \Lambda_{T_G} \otimes \Q \longto \Lambda_{T_H} \otimes \Q$,
  then the following diagram commutes:
  \begin{equation*} \xymatrix{
    \NS( \Bun_H) \ar[rr]^{(\iota_H)^{\NS, \varphi_* \delta}} \ar[d]^{\varphi_*} && \NS( \Bun_{T_H}) \ar[d]^{\varphi_*}\\
    \NS( \Bun_G) \ar[rr]^{(\iota_G)^{\NS, \delta}}                        && \NS( \Bun_{T_G}) 
  } \end{equation*}
\end{rem}

\subsection{The pullback to torus bundles} \label{sc2torus}
Let $\Ldet = \Ldet_n$ be determinant of cohomology line bundle \cite{M-K}
on $\Bun_{\GL_n}$, whose fibre at a vector bundle $E$ on $C$ is
$\det \H^*( E) = \det \H^0( E) \otimes \det \H^1( E)^{\dual}$.
\begin{lemma}
  Let $\xi$ be a line bundle of degree $d$ on $C$. Then the composition
  \begin{equation*}
    \Pic( \Bun_{\Gm}) \longto[ t_{\xi}^*] \Pic( \Bun_{\Gm}^0) \longto[ c_{\Gm}] \NS( \Bun_{\Gm}) = \Z \oplus \End_{J_C}
  \end{equation*}
  maps $\Ldet$ to $1 - g_C + d \in \Z$ and $-\id_{J_C} \in \End_{J_C}$.
\end{lemma}
\begin{proof}
  For any line bundle $L$ on $C$ and any point $p \in C( k)$, we have a canonical exact sequence
  \begin{equation*}
    0 \longto L( -p) \longto L \longto L_p \longto 0
  \end{equation*}
  of coherent sheaves on $C$. Varying $L$ and taking the determinant of cohomology, we see that the two line bundles
  $\Ldet$ and $t_{\O( -p)}^* \Ldet$ on $\Bun_{\Gm}^0$ have the same image in the second summand $\End J_C$ of $\NS( \Bun_{\Gm})$.
  Thus the image of $t_{\xi}^* \Ldet$ in $\End J_C$ does not depend on $\xi$; this image is $-\id_{J_C}$ because the
  principal polarisation $\phi_{ \Theta}: J_C \longto J_C^{\vee}$ is essentially given by the dual of the line bundle $\Ldet$.
  
  The weight of $t_{\xi}^* \Ldet$ at a line bundle $L$ of degree $0$ on $C$ is the Euler characteristic
  of $L \otimes \xi$, which is indeed $1 - g_C + d$ by Riemann--Roch theorem. 
\end{proof}
Let $\iota: T_{\SL_n} \hookrightarrow \SL_n$ be the inclusion of the maximal torus
$T_{\SL_n} := \Gm^n \cap \SL_n$, where $\Gm^n \subseteq \GL_n$ as diagonal matrices.
Then the cocharacter lattice $\Lambda_{T_{\SL_n}}$ is the group of all
$d = (d_1, \ldots, d_n) \in \Z^n$ with $d_1 + \cdots + d_n = 0$. The basic inner
product $b_{\SL_n}: \Lambda_{T_{\SL_n}} \otimes \Lambda_{T_{\SL_n}} \longto \Z$ is the
restriction of the standard scalar product on $\Z^n$.
\begin{cor} \label{cor:sl}
  Let $\xi$ be a principal $T_{\SL_n}$--bundle of degree
  $d \in \Lambda_{T_{\SL_n}}$ on $C$. Then the composition
  \begin{equation*}
    \Pic( \Bun_{\SL_n}) \longto[ \iota^*] \Pic( \Bun_{T_{\SL_n}}) \longto[ t_{\xi}^*] \Pic( \Bun^0_{T_{\SL_n}}) \xrightarrow{c_{T_{\SL_n}}} \NS( \Bun_{T_{\SL_n}})
  \end{equation*}
  maps $\Ldet$ to $b_{\SL_n}( d \otimes \_): \Lambda_{T_{\SL_n}} \longto \Z$ 
  and $-\id_{J_C} \cdot b_{\SL_n}: \Lambda_{T_{\SL_n}} \otimes \Lambda_{T_{\SL_n}} \longto \End J_C$.
\end{cor}
\begin{proof}
  Since the determinant of cohomology takes direct sums to tensor products, the pullback of $\Ldet_n$ to $\Bun_{\Gm^n}$ is isomorphic to
  $\pr_1^* \Ldet_1 \otimes \cdots \otimes \pr_n^* \Ldet_1$, where $\pr_{\nu}: \Gm^n \twoheadrightarrow \Gm$ is the projection onto
  the $\nu$th factor. Now use the previous lemma to compute the image of $\Ldet_n$ in $\NS( \Bun_{\Gm^n})$ and then restrict to $\NS( \Bun_{T_{\SL_n}})$.
\end{proof}
\begin{cor} \label{SL:dynkin}
  If $\rho: \SL_2 \longto \SL( V)$ has Dynkin index $d_{\rho}$, then the 
  pullback $\rho^*: \Pic( \Bun_{\SL( V)}) \longto \Pic( \Bun_{\SL_2})$ 
  maps $\Ldet$ to $(\Ldet_2)^{\otimes d_{\rho}}$.
\end{cor}
\begin{proof}
  Let $\iota: T_{\SL( V)} \hookrightarrow \SL( V)$ be the inclusion of a maximal torus that contains
  the image of the standard torus $T_{\SL_2} \hookrightarrow \SL_2$. The diagram
  \begin{equation*} \xymatrix{
    \Pic( \Bun_{\SL( V)}) \ar[r]^{\iota^*} \ar[d]^{\rho^*} & \Pic( \Bun_{T_{\SL( V)}}) \ar[r]^{t_{\rho_*( \xi)}^*} \ar[d]^{\rho^*}
      & \Pic( \Bun^0_{T_{\SL( V)}}) \ar[r]^{c_{T_{\SL( V)}}} \ar[d]^{\rho^*} & \NS( \Bun_{T_{\SL( V)}}) \ar[d]^{\rho^*}\\
    \Pic( \Bun_{\SL_2}) \ar[r]^{\iota^*} & \Pic( \Bun_{T_{\SL_2}}) \ar[r]^{t_{\xi}^*} & \Pic( \Bun^0_{T_{\SL_2}}) \ar[r]^{c_{T_{\SL_2}}} & \NS( \Bun_{T_{\SL_2}})
  } \end{equation*}
  commutes for each principal $T_{\SL_2}$--bundle $\xi$ on $C$.
  We choose $\xi$ in such a way that $\deg( \xi) \in \Lambda_{T_{\SL_2}} \cong \Z$ is nonzero if $g_C = 0$. Then the composition
  \begin{equation*}
    c_{T_{\SL_2}} \circ t_{\xi}^* \circ \iota^*: \Pic( \Bun_{\SL_2}) \longto \NS( \Bun_{T_{\SL_2}}) 
  \end{equation*}
  of the lower row is injective according to Theorem \ref{faltings} and 
  Corollary \ref{cor:sl}. The latter moreover implies that the two elements
  $\rho^*( \Ldet)$ and $(\Ldet_2)^{\otimes d_{\rho}}$ in $\Pic( \Bun_{\SL_2})$ have the same image in $\NS( \Bun_{T_{\SL_2}})$.
\end{proof}
Now suppose that the reductive group $G$ is simply connected and almost simple. We denote by $\O_{\Gr_G}( 1)$ the unique generator of $\Pic( \Gr_G)$
that is ample on every closed subscheme, and by $\O_{\Gr_G}( n)$ its $n$th tensor power for $n \in \Z$.

Over $k = \C$, the following is proved by a different method in section 5 of \cite{kumar-narasimhan-ramanathan}.
\begin{prop}[Kumar-Narasimhan-Ramanathan] \label{prop:dynkin}
  If $\rho: G \longto \SL( V)$ has Dynkin index $d_{\rho}$, then 
  $\rho^*: \Pic( \Gr_{\SL( V)}) \longto \Pic( \Gr_G)$ maps $\O( 1)$
  to $\O_{\Gr_G}( d_{\rho})$.
\end{prop}
\begin{proof}
  Let $\varphi: \SL_2 \longto G$ be given by a short coroot.
  Then  $d_{\varphi} = 1$ by definition, and \cite{faltings}
  implies that $\varphi^*: \Pic( \Gr_G) \longto \Pic( \Gr_{\SL_2})$
  maps $\O( 1)$ to $\O( 1)$, for example because
  $\varphi^*: \Pic( \Bun_G) \longto \Pic( \Bun_{\SL_2})$ preserves
  central charges according to their definition \cite[p. 59]{faltings}. 
  Hence it suffices to prove the claim for $\rho \circ \varphi$ 
  instead of $\rho$. This case follows from Corollary \ref{SL:dynkin},
  since $\glue_{p, z}^*( \Ldet_n) \cong \O_{\Gr_{\SL_n}}( -1)$.
\end{proof}
As in Subsection \ref{subsection:Grass}, we assume given an exact 
sequence of reductive groups
\begin{equation*}
  1 \longto G \longto \Ghat \longto[ \dt] \Gm \longto 1
\end{equation*}
with $G$ simply connected, and a line bundle $L$ on $C$.
\begin{cor} \label{charge_welldef}
  Suppose that $G$ is almost simple. Then the isomorphism
  \begin{equation*}
    (t_{\xi} \circ \glue_{p, z, \delta})^*: \relPic( \Bun_{\Ghat, L}) \longto[ \sim] \relPic( \Gr_G)
  \end{equation*}
  constructed in Subsection \ref{subsection:Grass} does not depend on 
  the choice of $p$, $z$, $\xi$ or $\delta$.
\end{cor}
We say that a line bundle on $\Bun_{\Ghat, L}$ has \emph{central charge $n \in \Z$}
if this isomorphism maps it to $\O_{\Gr_G}( n)$; this is consistent with the standard
central charge of line bundles on $\Bun_G$, as defined for example in \cite{faltings}.
\begin{proof}
  If $\rho: G \longto \SL( V)$ is a nontrivial representation, then $d_{\rho} > 0$,
  as explained in Remark \ref{simply_connected_products}(iii).
  Using Proposition \ref{prop:dynkin}, this implies that
  \begin{equation*}
    \rho^*: \Pic( \Gr_{\SL( V)}) \longto \Pic( \Gr_G)
  \end{equation*}
  is injective. Due to Remark \ref{glue_functorial},
  it thus suffices to check that
  \begin{equation*}
    \glue_{p, z}^*: \relPic( \Bun_{\SL( V)}) \longto[ \sim] \relPic( \Gr_{\SL( V)})
  \end{equation*}
  does not depend on $p$ or $z$. This is clear,
  since it maps $\Ldet$ to $\O_{\Gr_{\SL( V)}}( -1)$.
\end{proof}
The chosen maximal torus $\iota_G: T_G \hookrightarrow G$ induces maximal tori $\iota_{\Ghat}: T_{\Ghat} \hookrightarrow \Ghat$ and
$\iota_{\Gad}: T_{\Gad} \hookrightarrow \Gad$ compatible with the canonical maps $G \hookrightarrow \Ghat \twoheadrightarrow \Gad$.
Given a principal $T_{\Ghat}$--bundle $\xihat$ on $C$ and an isomorphism 
$\dt_* \xihat \cong L$, the composition
\begin{equation*}
  \Bun^0_{T_G} \longto[ t_{\xihat}] \Bun_{T_{\Ghat}} \longto[(\iota_{\Ghat})_*] \Bun_{\Ghat}
\end{equation*}
factors naturally through a $1$--morphism
\begin{equation} \label{iota_xi}
  \iota_{\xihat}: \Bun^0_{T_G} \longto \Bun_{\Ghat, L}.
\end{equation}
\begin{rem} \label{xihat_functorial}
  Given a representation $\rho: \Gad \longto \SL( V)$, let
  $\iota: T_{\SL( V)} \hookrightarrow \SL( V)$ be a maximal
  torus containing $\rho( T_{\Gad})$. Then the diagram
  \begin{equation*} \xymatrix{
    \Bun^0_{T_G} \ar[rr]^{\iota_{\xihat}} \ar[d]^{\rho_*} && \Bun_{\Ghat, L} \ar[d]^{\rho_*}\\
    \Bun^0_{T_{\SL( V)}} \ar[r]^{t_{\rho_* \xihat}} & \Bun_{T_{\SL( V)}} \ar[r]^{\iota_*} & \Bun_{\SL( V)}
  } \end{equation*}
  is $2$--commutative, by construction of $\iota_{\xihat}$.
\end{rem}
\begin{prop} \label{prop:twisted} \begin{itemize}
 \item[i)] $\Gamma( \Bun_{\Ghat, L}, \O_{\Bun_{\Ghat, L}}) = k$.
 \item[ii)] There is a canonical isomorphism
  \begin{equation*}
    c_G: \relPic( \Bun_{\Ghat, L}) \longto[ \sim] \NS( \Bun_G).
  \end{equation*}
 \item[iii)] For all choices of $\iota_G: T_G \hookrightarrow G$ and of $\xihat$, the diagram
  \begin{equation*} \xymatrix{
    \relPic( \Bun_{\Ghat, L}) \ar[rr]^{\iota_{\xihat}^*} \ar[d]^{c_G} && \relPic( \Bun^0_{T_G}) \ar[d]^{c_{T_G}}\\
    \NS( \Bun_G) \ar[rr]^{(\iota_G)^{\NS, \deltabar}} && \NS( \Bun_{T_G})
  } \end{equation*}
  commutes; here $\deltabar \in \Lambda_{T_{\Gad}}$ denotes the image  of $\deltahat := \deg \xihat \in \Lambda_{T_{\Ghat}}$.
\end{itemize} \end{prop}
\begin{proof}
  We start with the special case that $G$ is almost simple.
  Here part (i) of the proposition is just equation \eqref{fcts_simple} 
  from Subsection \ref{subsection:Grass}.

  We let $c_G$ send the line bundle of central charge $1$ to the
  basic inner product $b_G \in \NS( \Bun_G)$. Due to Theorem
  \ref{faltings}(i), Proposition \ref{prop:grass}, Corollary
  \ref{charge_welldef} and Remark \ref{simply_connected_products}(ii),
  this defines a canonical isomorphism, and hence proves (ii).

  To see that the diagram in (iii) then commutes, we choose a nontrivial 
  representation $\rho: \Gad \longto \SL( V)$. We note the functorialities,
  with respect to $\rho$, according to Remark \ref{xihat_functorial},
  Remark \ref{glue_functorial}, Proposition \ref{prop:dynkin},
  Remark \ref{iotaG_functorial} and Remark \ref{c_T:functorial}.
  In view of these, comparing Corollary \ref{cor:sl} and
  Definition \ref{def:iota_G} shows that the two images of
  $\rho^* \Ldet \in \Pic( \Bun_{\Ghat, L})$ in $\NS( \Bun_{T_G})$
  coincide. Since the former generates a subgroup of finite index
  and the latter is torsionfree, the diagram in (iii) commutes.

  For the general case, we use the unique decomposition
  \begin{equation*}
    G = G_1 \times \cdots \times G_r
  \end{equation*}
  into simply connected and almost simple factors $G_i$. As $\Ghat$ is generated by its center and $G$, every normal subgroup in $G$ is still normal in $\Ghat$.
  Let $\Ghat_i$ denote the quotient of $\Ghat$ modulo the closed normal subgroup $\prod_{j \neq i} G_j$; then
  \begin{equation*} \xymatrix{
    0 \ar[r] & G \ar[r] \ar@{->>}[d]^{\pr_i} & \Ghat \ar[r] \ar@{->>}[d] & \Gm \ar[r] \ar@{=}[d] & 0\\
    0 \ar[r] & G_i \ar[r] & \Ghat_i \ar[r]^{\dt_i} & \Gm \ar[r] & 0
  } \end{equation*}
  is a morphism of short exact sequences. Since the resulting diagram
  \begin{equation*} \xymatrix{
    \Ghat \ar[d]_{\dt} \ar[r] & \prod_i \Ghat_i \ar[d]^{\prod \dt_i}\\
    \Gm \ar[r]^{\diag} & \Gm^r
  } \end{equation*}
  is cartesian, it induces an equivalence of moduli stacks
  \begin{equation}\label{eqneu}
    \Bun_{\Ghat, L} \longto[ \sim] \Bun_{\Ghat_1, L} \times \cdots \times \Bun_{\Ghat_r, L}
  \end{equation}
  due to Lemma \ref{stacks:cartesian}. We note that equation
  \eqref{fcts_simple}, Lemma \ref{lemma:Pic}(i), Lemma \ref{see-saw},
  Remark \ref{simply_connected_products}(i) and  Corollary
  \ref{torus:product} ensure that various constructions are
  compatible with the products in \eqref{eqneu}. Therefore, the
  general case follows from the already treated almost simple case.
\end{proof}

\section{The reductive case} \label{section:reductive}
In this section, we finally describe the Picard functor $\relPic( \Bun_G^d)$ for
any reductive group $G$ over $k$ and any $d \in \pi_1( G)$. We denote
\begin{itemize}
 \item by $\zeta: Z^0 \hookrightarrow G$ the (reduced)
  identity component of the center $Z \subseteq G$, and
 \item by $\pi: \Gtilde \longto G$ the universal cover
  of $G' := [G\, , G] \subseteq G$.
\end{itemize}
Our strategy is to descend along the central isogeny
\begin{equation*}
  \zeta \cdot \pi: Z^0 \times \Gtilde \longto G,
\end{equation*}
applying the previous two sections to $Z^0$ and to $\Gtilde$, 
respectively. The $1$--morphism of moduli stacks given by such a central isogeny is a torsor
under a group stack; Subsection \ref{stacktorsors} explains descent of 
line bundles along such torsors, generalising the method introduced by Laszlo
\cite{laszlo} for quotients of $\SL_n$. In Subsection 
\ref{subsection:NS_reductive}, we define combinatorially what will be the discrete torsionfree part
of $\relPic( \Bun_G^d)$; finally, these Picard functors and their 
functoriality in $G$ are described in Subsection \ref{mainresult}.

The following notation is used throughout this section. The reductive group $G$ yields semisimple groups and central isogenies
\begin{equation*}
  \Gtilde \twoheadrightarrow G' \twoheadrightarrow \Gbar := G/Z^0 \twoheadrightarrow \Gad := G/Z.
\end{equation*}
We denote by $\dbar \in \pi_1( \Gbar) \subseteq \pi_1( \Gad)$ the image of $d \in \pi_1( G)$.
The choice of a maximal torus $\iota_G: T_G \hookrightarrow G$ induces maximal tori and isogenies
\begin{equation*}
  T_{\Gtilde} \twoheadrightarrow T_{G'} \twoheadrightarrow T_{\Gbar} \twoheadrightarrow T_{\Gad}.
\end{equation*}
Their cocharacter lattices are hence subgroups of finite index
\begin{equation*}
  \Lambda_{T_{\Gtilde}} \hookrightarrow \Lambda_{T_{G'}} \hookrightarrow \Lambda_{T_{\Gbar}} \hookrightarrow \Lambda_{T_{\Gad}}.
\end{equation*}
The central isogeny $\zeta \cdot \pi$ makes $\Lambda_{Z^0} \oplus \Lambda_{T_{\Gtilde}}$ a subgroup of finite index in $\Lambda_{T_G}$.

\subsection{Torsors under a group stack} \label{stacktorsors}
All stacks in this subsection are stacks over $k$, and all morphisms are over $k$.
Following \cite{breen, laszlo}, we recall the notion of a torsor under a group stack.

Let $\calG$ be a group stack. We denote by $1$ the unit object in $\calG$,
and by $g_1 \cdot g_2$ the image of two objects $g_1$ and $g_2$ under
the multiplication $1$--morphism $\calG \times \calG \longto \calG$.
\begin{defn}
  An \emph{action} of $\calG$ on a $1$--morphism of stacks
  $\Phi: \calX \longto \calY$ consists of a $1$--morphism
  \begin{equation*}
    \calG \times \calX \longto \calX, \qquad (g, x) \longmapsto g \cdot x,
  \end{equation*}
  and of three $2$--morphisms, which assign
  to each $k$--scheme $S$ and each object
  \begin{align*}
    x & \text{ in } \calX( S)
      & \text{an isomorphism } 1 \cdot x
      & \longto[ \sim] x \text{ in } \calX( S),\\
    (g, x) & \text{ in } (\calG \times \calX)( S)
      & \text{an isomorphism } \Phi( g \cdot x)
      & \longto[ \sim] \Phi( x) \text{ in } \calY( S),\\
    (g_1,g_2,x) & \text{ in } (\calG \times \calG \times \calX)(S)
      & \text{an isomorphism } (g_1 \cdot g_2)\cdot x
      & \longto[ \sim] g_1 \cdot (g_2 \cdot x) \text{ in } \calX( S).
  \end{align*}
  These morphisms are required to satisfy the following five compatibility conditions: the two resulting isomorphisms
  \begin{align*}
    (g \cdot 1) \cdot x               & \longto[ \sim] g \cdot x \text{ in } \calX( S),\\
    (1 \cdot g) \cdot x               & \longto[ \sim] g \cdot x \text{ in } \calX( S),\\
    \Phi( 1 \cdot x)                  & \longto[ \sim] \Phi( x) \text{ in } \calY( S),\\
    \Phi( (g_1 \cdot g_2) \cdot x)    & \longto[ \sim] \Phi( x) \text{ in } \calY( S),\\ \text{and }
    (g_1 \cdot g_2 \cdot g_3) \cdot x & \longto[ \sim] g_1 \cdot ( g_2 \cdot ( g_3 \cdot x)) \text{ in } \calX( S), 
  \end{align*}
  coincide for all $k$--schemes $S$ and all objects
  $g, g_1, g_2, g_3$ in $\calG( S)$ and $x$ in $\calX( S)$.
\end{defn}

\begin{example} \label{action}
  Let $\varphi: G \longto H$ be a homomorphism of linear algebraic 
  groups over $k$, and let $Z$ be a closed subgroup in the center
  of $G$ with $Z \subseteq \ker( \varphi)$.
 
  Then the group stack $\Bun_Z$ acts on the $1$--morphism
  $\varphi_*: \Bun_G \longto \Bun_H$ via the tensor product
  $\_ \otimes \_: \Bun_Z \times \Bun_G \longto \Bun_G$.
\end{example}

From now on, we assume that the group stack $\calG$ is algebraic.
\begin{defn}
  A \emph{$\calG$--torsor} is a faithfully flat $1$--morphism of 
  algebraic stacks $\Phi: \calX \longto \calY$ together with
  an action of $\calG$ on $\Phi$ such that the resulting $1$--morphism
  \begin{equation*}
    \calG \times \calX \longto \calX \times_{\calY} \calX,
      \quad (g, x) \longmapsto (g \cdot x, x)
  \end{equation*}
  is an isomorphism.
\end{defn}

\begin{example} \label{isogeny2torsor} \upshape
  Suppose that $\varphi: G \twoheadrightarrow H$ is a central isogeny of reductive groups with kernel $\mu$.
  For each $d \in \pi_1( G)$, the $1$--morphism
  \begin{equation} \label{torsor?}
    \varphi_*: \Bun_G^d \longto \Bun_H^e \qquad\qquad\qquad e := \varphi_*( d) \in \pi_1( H)
  \end{equation}
  is a torsor under the group stack $\Bun_{\mu}$, for the action described in example \ref{action}.
\end{example}

\begin{proof}
  The $1$--morphism $\varphi_*$ is faithfully flat by Lemma \ref{flat}. 
  The $1$--morphism
  \begin{equation*}
    \Bun_{\mu} \times \Bun_G \longto \Bun_G \times_{\Bun_H} \Bun_G, 
      \quad (L, E) \longmapsto (L \otimes E, E)
  \end{equation*}
  is an isomorphism due to Lemma \ref{stacks:cartesian}. Since 
  $\varphi_*: \pi_1( G) \longto \pi_1( H)$ is injective,
  $\Bun_G^d \subseteq \Bun_G$ is the inverse image of
  $\Bun_H^e \subseteq \Bun_H$ under $\varphi_*$; hence the restriction
  \begin{equation*}
    \Bun_{\mu} \times \Bun_G^d \longto \Bun_G^d \times_{\Bun_H^e} \Bun_G^d
  \end{equation*}
  is an isomorphism as well.
\end{proof}

\begin{defn}
  Let $\Phi_{\nu}: \calX_{\nu} \longto \calY_{\nu}$ be a $\calG$--torsor 
  for $\nu = 1, 2$. A \emph{morphism of $\calG$--torsors} from $\Phi_1$
  to $\Phi_2$ consists of two $1$--morphisms
  \begin{equation*}
    A: \calX_1 \longto \calX_2 \qquad\text{and}\qquad B: \calY_1 \longto \calY_2
  \end{equation*}
  and of two $2$--morphisms, which assign to each $k$--scheme $S$
  and each object
  \begin{align*}
    x & \text{ in } \calX_1( S)
      & \text{an isomorphism } \Phi_2 A( x)
      & \longto[ \sim] B \Phi_1( x) \text{ in } \calY_2( S),\\
    (g, x) & \text{ in } (\calG \times \calX_1)( S)
      & \text{an isomorphism } A( g \cdot x)
      & \longto[ \sim] g \cdot A( x) \text{ in } \calX_2( S).
  \end{align*}
  These morphisms are required to satisfy the following three compatibility conditions: the two resulting isomorphisms
  \begin{align*}
    A( 1 \cdot x) & \longto[ \sim] A( x) \text{ in } \calX_2( S),\\
    \Phi_2 A( g \cdot x) & \longto[ \sim] B \Phi_1( x) \text{ in } \calY_2( S)\\ \text{and }
    A( (g_1 \cdot g_2) \cdot x) & \longto[ \sim] g_1 \cdot ( g_2 \cdot A( x)) \text{ in } \calX_2( S)
  \end{align*}
  coincide for all $k$--schemes $S$ and all objects
  $g, g_1, g_2$ in $\calG( S)$ and $x$ in $\calX_1( S)$.
\end{defn}

\begin{example} \label{pullback2torsors}
  Let a cartesian square of reductive groups over $k$
  \begin{equation*} \xymatrix{
    G_1 \ar[r]^{\alpha} \ar[d]_{\varphi_1} & G_2 \ar[d]^{\varphi_2}\\
    H_1 \ar[r]^{\beta} & H_2
  } \end{equation*}
  be given. Suppose that $\varphi_1$ and $\varphi_2$ are central isogenies, and denote their common kernel by $\mu$. For each $d_1 \in \pi_1( G_1)$, the diagram
  \begin{equation*} \xymatrix{
    \Bun_{G_1}^{d_1} \ar[r]^{\alpha_*} \ar[d]_{(\varphi_1)_*} & \Bun_{G_2}^{d_2} \ar[d]^{(\varphi_2)_*} && d_2 := \alpha_*( d_1) \in \pi_1( G_2)\\
    \Bun_{H_1}^{e_1} \ar[r]^{\beta_*} & \Bun_{H_2}^{e_2} && e_{\nu} := (\varphi_{\nu})_*( d_{\nu}) \in \pi_1( H_{\nu})
  } \end{equation*}
  is then a morphism of torsors under the group stack $\Bun_{\mu}$.
\end{example}

\begin{prop} \label{torsors:cartesian}
  Let a $\calG$--torsor $\Phi_{\nu}: \calX_{\nu} \longto \calY_{\nu}$
  with $\Gamma( \calX_{\nu}, \O_{\calX_{\nu}}) = k$ be given for
  $\nu = 1, 2$, together with a morphism of $\calG$--torsors
  \begin{equation*} \xymatrix{
    \calX_1 \ar[r]^A \ar[d]_{\Phi_1} & \calX_2 \ar[d]^{\Phi_2}\\
    \calY_1 \ar[r]^B & \calY_2
  } \end{equation*}
  such that the induced morphism of Picard functors
  $A^*: \relPic( \calX_2) \longto \relPic( \calX_1)$
  is injective. Then the diagram of Picard functors
  \begin{equation*} \xymatrix{
    \relPic( \calX_1) & \relPic( \calX_2) \ar[l]_-{A^*}\\
    \relPic( \calY_1) \ar[u]^{\Phi_1^*} & \relPic( \calY_2) \ar[u]_{\Phi_2^*} \ar[l]_-{B^*}
  } \end{equation*}
  is a pullback square.
\end{prop}
\begin{proof}
  The proof of \cite[Theorem 5.7]{laszlo} generalises to this situation as follows.

  Let $S$ be a scheme of finite type over $k$. For a line bundle $\calL$
  on $S \times \calX_{\nu}$, we denote by $\Lin^{\calG}( \calL)$ the set
  of its $\calG$--linearisations, cf. \cite[Definition 2.8]{laszlo}.
  According to Lemma \ref{lemma:Pic}(i), each automorphism of $\calL$
  comes from $\Gamma( S, \O_S^*)$ and hence respects each linearisation
  of $\calL$. Thus \cite[Theorem 4.1]{laszlo} provides a canonical
  bijection between the set $\Lin^{\calG}( \calL)$ and the fibre of
  \begin{equation*}
    \Phi_{\nu}^*: \Pic( S \times \calY_{\nu}) \longto \Pic( S \times \calX_{\nu})
  \end{equation*}
  over the isomorphism class of $\calL$.

  Let $\calT$ be an algebraic stack over $k$. We denote for the moment by
  $\calPic( \calT)$ the groupoid of line bundles on $\calT$ and their isomorphisms.
  Lemma \ref{lemma:Pic}(i) and Corollary \ref{Pic:injective} show that the functor
  \begin{equation*}
    A^*: \calPic( \calT \times \calX_2) \longto \calPic( \calT \times \calX_1)
  \end{equation*}
  is fully faithful for every $\calT$. We recall that an element in $\Lin^{\calG}( \calL)$ is an isomorphism
  in $\calPic( \calG \times S \times \calX_{\nu})$ between two pullbacks of $\calL$ such that certain induced diagrams
  in $\calPic( S \times \calX_{\nu})$ and in $\calPic( \calG \times \calG \times S \times \calX_{\nu})$ commute.
  Thus it follows for all $\calL \in \Pic( S \times \calX_2)$ that the canonical map
  \begin{equation*}
    A^*: \Lin^{\calG}( \calL) \longto \Lin^{\calG}( A^* \calL)
  \end{equation*}
  is bijective. Hence the diagram of abelian groups
  \begin{equation*} \xymatrix{
    \Pic( S \times \calX_1) & \Pic( S \times \calX_2) \ar[l]_-{A^*}\\
    \Pic( S \times \calY_1) \ar[u]^{\Phi_1^*} & \Pic( S \times \calY_2) \ar[u]_{\Phi_2^*} \ar[l]_-{B^*}
  } \end{equation*}
  is a pullback square, as required.
\end{proof}

\subsection{N\'{e}ron--Severi groups $\NS( \Bun_G^d)$ for reductive $G$} 
\label{subsection:NS_reductive}
\begin{defn} \label{NS:reductive}
  The \emph{N\'{e}ron--Severi group} $\NS( \Bun_G^d)$ is the subgroup
  \begin{equation*}
    \NS( \Bun_G^d) \subseteq \NS( \Bun_{Z^0}) \oplus \NS( \Bun_{\Gtilde})
  \end{equation*}
  of all triples $l_Z: \Lambda_{Z^0} \longto \Z$,
  $b_Z: \Lambda_{Z^0} \otimes \Lambda_{Z^0} \longto \End J_C$ and
  $b: \Lambda_{T_{\Gtilde}} \otimes \Lambda_{T_{\Gtilde}} \longto \Z$
  with the following properties:
  \begin{enumerate}
   \item \label{extend_weights} For every lift $\deltabar \in \Lambda_{T_{\Gbar}}$ of $\dbar \in \pi_1( \Gbar)$, the direct sum
    \begin{equation*}
      l_Z \oplus b( -\deltabar \otimes \_): \Lambda_{Z^0} \oplus \Lambda_{T_{\Gtilde}} \longto \Z
    \end{equation*}
    is integral on $\Lambda_{T_G}$.
   \item \label{extend_form} The orthogonal direct sum
    \begin{equation*}
      b_Z \perp (\id_{J_C} \cdot b): (\Lambda_{Z^0} \oplus \Lambda_{T_{\Gtilde}}) \otimes (\Lambda_{Z^0} \oplus \Lambda_{T_{\Gtilde}}) \longto \End J_C
    \end{equation*}
    is integral on $\Lambda_{T_G} \otimes \Lambda_{T_G}$.
  \end{enumerate}
\end{defn}
\begin{lemma} \label{NS:welldef}
  If condition \ref{extend_weights} above holds for one lift $\deltabar \in \Lambda_{T_{\Gbar}}$ of $\dbar \in \pi_1( \Gbar)$,
  then it holds for every lift $\deltabar \in \Lambda_{T_{\Gbar}}$ of the same element $\dbar \in \pi_1( \Gbar)$.
\end{lemma}
\begin{proof}
  Any two lifts $\deltabar$ of $\dbar$ differ by some element
  $\lambda \in \Lambda_{T_{\Gtilde}}$. Lemma \ref{extend_b} states
  in particular that
  \begin{equation*}
    b( -\lambda \otimes \_): \Lambda_{T_{\Gtilde}} \longto \Z
  \end{equation*}
  is integral on $\Lambda_{T_{\Gbar}}$, and hence admits an extension
  $\Lambda_{T_G} \longto \Z$ that vanishes on $\Lambda_{Z^0}$.
\end{proof}
\begin{rem}
  If $G$ is simply connected, then $\NS( \Bun_G^0)$ coincides with the group $\NS( \Bun_G)$ of definition \ref{NS:simply_connected}.
  If $G = T$ is a torus, then $\NS( \Bun_T^d)$ coincides for all $d \in \pi_1( T)$ with the group $\NS( \Bun_T)$ of definition \ref{NS:torus}.
\end{rem}
\begin{rem}
  The Weyl group $W$ of $(G, T_G)$ acts trivially on $\NS( \Bun_G^d)$.
  Hence the group $\NS( \Bun_G^d)$ does not depend on the choice of 
  $T_G$; cf. Subsection \ref{sc:NS}.
\end{rem}
\begin{defn} \label{def:iota_G_reductive}
  Given a lift $\delta \in \Lambda_{T_G}$ of $d \in \pi_1( G)$, the homomorphism
  \begin{equation*}
    (\iota_G)^{\NS, \delta}: \NS( \Bun_G^d) \longto \NS( \Bun_{T_G})
  \end{equation*}
  sends $(l_Z, b_Z) \in \NS( \Bun_{Z^0})$ and $b \in \NS( \Bun_{\Gtilde})$ to the pair
  \begin{equation*}
    l_Z \oplus b( -\deltabar \otimes \_): \Lambda_G \longto \Z \quad\text{and}\quad
    b_Z \perp (\id_{J_C} \cdot b): \Lambda_{T_G} \otimes \Lambda_{T_G} \longto \End J_C
  \end{equation*}
  where $\deltabar \in \Lambda_{T_{\Gbar}}$ denotes the image of $\delta$.
\end{defn}
Note that this definition agrees with the earlier definition \ref{def:iota_G} in the cases covered by both,
namely $G$ simply connected and $\delta \in \Lambda_{T_G}$.
\begin{lemma} \label{NS:cartesian}
  Given a lift $\delta \in \Lambda_{T_G}$ of $d \in \pi_1( G)$, the diagram
  \begin{equation*} \xymatrix{
    \NS( \Bun_G^d) \ar[rrr]^{(\iota_G)^{\NS, \delta}} \ar@{^{(}->}[d] &&& \NS( \Bun_{T_G}) \ar[d]^{(\zeta \cdot \pi)^*}\\
    \NS( \Bun_{Z^0}) \oplus \NS( \Bun_{\Gtilde}) \ar[rr]^-{\id \oplus (\iota_{\Gtilde})^{\NS, \deltabar}}
      && \NS( \Bun_{Z^0}) \oplus \NS( \Bun_{T_{\Gtilde}}) \ar@{^{(}->}[r] & \NS( \Bun_{Z^0 \times T_{\Gtilde}})
  } \end{equation*}
  is a pullback square; here $\deltabar \in \Lambda_{T_{\Gad}}$ again denotes the image of $\delta$.
\end{lemma}
\begin{proof}
  This follows directly from the definitions.
\end{proof}
Let $e \in \pi_1( H)$ be the image of $d \in \pi_1( G)$ under a 
homomorphism of reductive groups $\varphi: G \longto H$. $\varphi$ 
induces a map $\varphi: \Gtilde \longto \Htilde$ between the
universal covers of their commutator subgroups. If $\varphi$
maps the identity component $Z^0_G$ in the center $Z_G$ of $G$
to the center $Z_H$ of $H$, then it induces an obvious pullback map
\begin{equation*}
  \varphi^*: \NS( \Bun_H^e) \longto \NS( \Bun_G^d)
\end{equation*}
which sends $l_Z$, $b_Z$ and $b$ simply to $\varphi^* l_Z$, $\varphi^* b_Z$ and $\varphi^* b$. This is a special case of the following
map, which $\varphi$ induces even without the hypothesis on the centers, and which also generalises the previous definition \ref{def:iota_G_reductive}.
\begin{defn} \label{NS_reductive:functorial}
  Choose a maximal torus $\iota_H: T_H \hookrightarrow H$ containing $\varphi( T_G)$, and a lift $\delta \in \Lambda_{T_G}$ of $d \in \pi_1( G)$;
  let $\eta \in \Lambda_{T_H}$ be the image of $\delta$. Then the map
  \begin{equation*}
    \varphi^{\NS, d}: \NS( \Bun_H^e) \longto \NS( \Bun_G^d)
  \end{equation*}
  sends $(l_Z, b_Z) \in \NS( \Bun_{Z^0_H})$ and $b \in \NS( \Bun_{\Htilde})$
  to the pullback along $\varphi: Z_G^0 \longto T_H$ of
  $(\iota_H)^{\NS, \eta}( l_Z, b_Z, b) \in \NS( \Bun_{T_H})$,
  together with $\varphi^* b \in \NS( \Bun_{\Gtilde})$.
\end{defn}
\begin{lemma}
  The map $\varphi^{\NS, d}$ does not depend
  on the choice of $T_G$, $T_H$ or $\delta$.
\end{lemma}
\begin{proof}
  Let $W_G$ denote the Weyl group of $(G, T_G)$. It acts trivially on $\Lambda_{Z^0_G}$, and without nontrivial coinvariants on $\Lambda_{T_{\Gtilde}}$;
  these two observations imply
  \begin{equation} \label{orthogonal}
    \Hom( \Lambda_{T_{\Gtilde}} \otimes \Lambda_{Z^0_G}, \Z)^{W_G} = 0.
  \end{equation}
  Lemma \ref{extend_b} states that $b$ is integral on
  $\Lambda_{T_{\Htilde}} \otimes \Lambda_{T_{\Hbar}}$; its
  composition with the canonical projection
  $\Lambda_{T_H} \twoheadrightarrow \Lambda_{T_{\Hbar}}$ is a 
  Weyl--invariant map $b_r: \Lambda_{T_{\Htilde}} \otimes \Lambda_{T_H} \longto \Z$.
  As explained in Subsection \ref{sc:NS}, Lemma \ref{from_Weyl} implies 
  that $\varphi^* b_r: \Lambda_{T_{\Gtilde}} \otimes \Lambda_{T_G} \longto \Z$
  is still Weyl--invariant; hence it vanishes on
  $\Lambda_{T_{\Gtilde}} \otimes \Lambda_{Z^0_G}$ by \eqref{orthogonal}.

  Any two lifts $\delta$ of $d$ differ by some element
  $\lambda \in \Lambda_{T_{\Gtilde}}$; then the two images of
  $( l_Z, b_Z, b) \in \NS( \Bun_H^e)$ in $\NS( \Bun_{T_H})$
  differ, according to the proof of Lemma \ref{NS:welldef},
  only by $b_r( -\lambda \otimes \_): \Lambda_{T_H} \longto \Z$.
  Thus their compositions with $\varphi: \Lambda_{Z^0_G} \longto \Lambda_{T_H}$
  coincide by the previous paragraph. This shows that the two images of
  $( l_Z, b_Z, b)$ have the same component in the direct summand
  $\Hom( \Lambda_{Z^0_G}, \Z)$ of $\NS( \Bun_G^d)$; since the other two 
  components do not involve $\delta$ at all, the independence on $\delta$ follows.

  The independence on $T_G$ and $T_H$ is then a consequence of Lemma 
  \ref{from_Weyl}, since the Weyl groups $W_G$ and $W_H$ act trivially
  on $\NS( \Bun_G^d)$ and on $\NS( \Bun_H^e)$.
\end{proof}
\begin{lemma} \label{iota^NS_functorial}
  For all maximal tori $\iota_G: T_G \hookrightarrow G$ and $\iota_H: T_H \hookrightarrow H$ with $\varphi( T_G) \subseteq T_H$,
  and all lifts $\delta \in \Lambda_{T_G}$ of $d \in \pi_1( G)$, the diagram
  \begin{equation*} \xymatrix{
    \NS( \Bun_H^e) \ar[rr]^{(\iota_H)^{\NS, \eta}} \ar[d]^{\varphi^{\NS, d}} && \NS( \Bun_{T_H}) \ar[d]^{\varphi^*}\\
    \NS( \Bun_G^d) \ar[rr]^{(\iota_G)^{\NS, \delta}} && \NS( \Bun_{T_G})
  } \end{equation*}
  commutes, with $\eta := \varphi_* \delta \in \Lambda_{T_H}$ and $e := \varphi_* d \in \pi_1( H)$ as in definition \ref{NS_reductive:functorial}.
\end{lemma}
\begin{proof}
  Given an element in $\NS( \Bun^e_H)$, we have to compare its two images in $\NS( \Bun_{T_G})$. The definition \ref{NS_reductive:functorial} of $\varphi^{\NS, d}$
  directly implies that both have the same pullback to $\NS( \Bun_{Z^0_G})$ and to $\NS( \Bun_{T_{\Gtilde}})$. Moreover, their components in the direct summand
  $\Hom^s( \Lambda_{T_G} \otimes \Lambda_{T_G}, \End J_C)$ of $\NS( \Bun_{T_G})$ are both Weyl--invariant due to Lemma \ref{from_Weyl}; thus equation
  \eqref{orthogonal} above shows that these components vanish on $\Lambda_{T_{\Gtilde}} \otimes \Lambda_{Z^0_G}$ and on $\Lambda_{Z^0_G} \otimes \Lambda_{T_{\Gtilde}}$. 
  Hence two images in question even have the same pullback to $\NS( \Bun_{Z^0_G \times T_{\Gtilde}})$. But $\Lambda_{Z^0_G} \oplus \Lambda_{T_{\Gtilde}}$
  has finite index in $\Lambda_{T_G}$.
\end{proof}
\begin{cor} \label{phi^NS_functorial}
  Let $\psi: H \longto K$ be another homomorphism of reductive groups, 
  and put $f := \psi_* e \in \pi_1( K)$. Then
  \begin{equation*}
    \varphi^{\NS, d} \circ \psi^{\NS, e} = (\psi \circ \varphi)^{\NS, d}: \NS( \Bun_K^f) \longto \NS( \Bun_G^d).
  \end{equation*}
\end{cor}
\begin{proof}
  According to the previous lemma, this equality holds after composition 
  with $(\iota_G)^{\NS, \delta}: \NS( \Bun_G^d) \longto \NS( \Bun_{T_G})$
  for any lift $\delta \in \Lambda_{T_G}$ of $d$. Due to the Lemma 
  \ref{injective} and Lemma \ref{NS:cartesian}, there is a lift $\delta$ 
  of $d$ such that $(\iota_G)^{\NS, \delta}$ is injective.
\end{proof}
We conclude this subsection with a more explicit description of $\NS( \Bun_G^d)$. It turns out that genus $g_C = 0$ is special.
This generalises the description obtained for $k = \C$ and $G$ semisimple by different methods in \cite[Section V]{teleman}.
\begin{prop} \label{NS_extension}
  Let $q: G \twoheadrightarrow G/G' =: G^{\ab}$ denote the maximal abelian quotient of $G$. Then the sequence of abelian groups
  \begin{equation*}
    0 \longto \NS( \Bun_{G^{\ab}}) \longto[ q^*] \NS( \Bun_G^d) \longto[ \pr_2] \NS( \Bun_{\Gtilde})
  \end{equation*}
  is exact, and the image of the map $\pr_2$ in it consists of
  all forms $b: \Lambda_{T_{\Gtilde}} \otimes \Lambda_{T_{\Gtilde}} \longto \Z$
  in $\NS( \Bun_{\Gtilde})$ that are integral
  \begin{itemize}
   \item on $\Lambda_{T_{\Gbar}} \otimes \Lambda_{T_{G'}}$, if $g_C \geq 1$;
   \item on $(\Z \deltabar) \otimes \Lambda_{T_{G'}}$ for a lift $\deltabar \in \Lambda_{T_{\Gbar}}$ of $\dbar \in \pi_1( \Gbar)$, if $g_C = 0$.
  \end{itemize}
  The condition does not depend on the choice of
  this lift $\deltabar$, due to Lemma \ref{extend_b}.
\end{prop}
\begin{proof}
  Since $q: Z^0 \longto G^{\ab}$ is an isogeny, $q^*$ is injective;
  it clearly maps into the kernel of $\pr_2$. Conversely, let
  $(l_Z, b_Z, b) \in \NS( \Bun_G^d)$ be in the kernel of $\pr_2$;
  this means $b = 0$. Then condition \ref{extend_weights} in the
  definition \ref{NS:reductive} of $\NS( \Bun_G^d)$ provides a map
  \begin{equation*} 
    l_Z \oplus 0: \Lambda_{T_G} \longto \Z
  \end{equation*}
  which vanishes on $\Lambda_{T_{\Gtilde}}$, and hence also on $\Lambda_{T_{G'}}$;
  thus it is induced from a map on $\Lambda_{T_G}/\Lambda_{T_{G'}} = \Lambda_{G^{\ab}}$.
  Similarly, condition \ref{extend_form} in the same definition provides a map
  $b_Z \perp 0$ on $\Lambda_{T_G} \otimes \Lambda_{T_G}$ which vanishes on
  $\Lambda_{T_{\Gtilde}} \otimes \Lambda_{T_G} + \Lambda_{T_G} \otimes \Lambda_{T_{\Gtilde}}$, and hence
  also on $\Lambda_{T_{G'}} \otimes \Lambda_{T_G} + \Lambda_{T_G} \otimes \Lambda_{T_{G'}}$;
  thus it is induced from a map on the quotient $\Lambda_{G^{\ab}} \otimes \Lambda_{G^{\ab}}$. This proves the exactness.

  Now let $b \in \NS( \Bun_G^d)$ be in the image of $\pr_2$. Then $b$ is integral on
  $(\Z \deltabar) \otimes \Lambda_{G'}$ by condition \ref{extend_weights}
  in definition \ref{NS:reductive}. If $g_C \geq 1$, then
  \begin{equation*}
    \_ \cdot \id_{J_C}: \Z \longto \End J_C
  \end{equation*}
  is injective with torsionfree cokernel; thus condition
  \ref{extend_form} in definition \ref{NS:reductive} implies that
  \begin{equation*}
    0 \oplus b: ( \Lambda_{Z^0} \oplus \Lambda_{T_{\Gtilde}}) \otimes \Lambda_{T_{\Gtilde}} \longto \Z
  \end{equation*}
  is integral on $\Lambda_{T_G} \otimes \Lambda_{T_{G'}}$ and hence, vanishing on $\Lambda_{Z^0} \subseteq \Lambda_{T_G}$, comes from a map on the quotient
  $\Lambda_{T_{\Gbar}} \otimes \Lambda_{T_{G'}}$. This shows that $b$ satisfies the stated condition.

  Conversely, suppose that $b \in \NS( \Bun_G^d)$ satisfies the stated condition.
  Then $b$ is integral on $(\Z \deltabar) \otimes \Lambda_{T_{G'}}$; since
  $\Lambda_{T_{G'}} \subseteq \Lambda_{T_G}$ is a direct summand,
  \begin{equation*}
    b( -\deltabar \otimes \_): \Lambda_{T_{G'}} \longto \Z
  \end{equation*}
  can thus be extended to $\Lambda_{T_G}$. We restrict it to a map
  $l_Z: \Lambda_{Z^0} \longto \Z$. In the case $g_C = 0$, the triple
  $(l_Z, 0, b)$ is in $\NS( \Bun_G^d)$ and hence an inverse image of $b$.

  It remains to consider $g_C \geq 1$. Then $b$ is by assumption integral
  on $\Lambda_{T_{\Gbar}} \otimes \Lambda_{T_{G'}}$, so composing it with the canonical
  subjection $\Lambda_{T_G} \twoheadrightarrow \Lambda_{T_{\Gbar}}$ defines a
  linear map $\Lambda_{T_G} \otimes \Lambda_{T_{G'}} \longto \Z$.
  Since $b$ is symmetric, this extends canonically to a symmetric linear map from
  \begin{equation*}
    \Lambda_{T_G} \otimes \Lambda_{T_{G'}} + \Lambda_{T_{G'}} \otimes \Lambda_{T_G} \subseteq \Lambda_{T_G} \otimes \Lambda_{T_G}
  \end{equation*}
  to $\Z$. It can be extended further to a symmetric linear map from
  $\Lambda_{T_G} \otimes \Lambda_{T_G}$ to $\Z$, because $\Lambda_{T_{G'}} \subseteq \Lambda_{T_G}$
  is a direct summand. Multiplying it with $\id_{J_C}$ and restricting to $\Lambda_{Z^0}$
  defines an element $b_Z \in \Hom^s( \Lambda_{Z^0} \otimes \Lambda_{Z^0}, \End J_C)$.
  By construction, the triple $(l_Z, b_Z, b)$ is in $\NS( \Bun_G^d)$ and hence an inverse image of $b$.
\end{proof}
In particular, the free abelian group $\NS( \Bun_G^d)$ has rank
\begin{equation*}
  \rank \NS( \Bun_G^d) = r + r \cdot \rank \NS( J_C) + \frac{r(r-1)}{2} \cdot \rank \End( J_C) + s
\end{equation*}
if $G^{\ab} \cong \Gm^r$ is a torus of rank $r$, and $\Gad$ contains $s$ simple factors.

\subsection{Proof of the main result} \label{mainresult}
\begin{thm} \label{mainthm}
  \begin{itemize}  
   \item[i)] $\Gamma( \Bun_G^d, \O_{\Bun_G^d}) = k$.
   \item[ii)] The functor $\relPic( \Bun_G^d)$ is representable by a 
    $k$--scheme locally of finite type.
   \item[iii)] There is a canonical exact sequence
    \begin{equation*}
      0 \longto \relHom( \pi_1( G), J_C) \longto[ j_G] \relPic( \Bun_G^d) \longto[ c_G] \NS( \Bun_G^d) \longto 0
    \end{equation*}
    of commutative group schemes over $k$.
   \item[iv)] For every homomorphism of reductive groups
    $\varphi: G \longto H$, the diagram
    \begin{equation*} \xymatrix{
      0\ar[r] &\relHom(\pi_1(H),J_C) \ar[r]^-{j_H}\ar[d]^{\varphi^*} &\relPic(\Bun_H^e) \ar[r]^{c_H}\ar[d]^{\varphi^*} &\NS(\Bun_H^e) \ar[r]\ar[d]^{\varphi^{\NS, d}}&0\\
      0\ar[r] &\relHom(\pi_1(G),J_C) \ar[r]^-{j_G}                  &\relPic(\Bun_G^d) \ar[r]^{c_G}                 &\NS(\Bun_G^d) \ar[r]                     &0
    } \end{equation*}
    commutes; here $e := \varphi_*( d) \in \pi_1( H)$.
  \end{itemize}
\end{thm}
\begin{proof}
  We record for later use the commutative square of abelian groups
  \begin{equation*} \xymatrix{
    \pi_1( G) && \Lambda_{T_G} \ar@{->>}[ll]_{\pr}\\
    \Lambda_{Z^0} \ar[u]_{\zeta_*} && \Lambda_{Z^0 \times T_{\Gtilde}}. \ar[u]_{(\zeta \cdot \pi)_*} \ar@{->>}[ll]_{\pr_1}
  } \end{equation*}
  The mapping cone of this commutative square
  \begin{equation} \label{lattices:exact}
    0 \longto \Lambda_{Z^0} \oplus \Lambda_{T_{\Gtilde}} \longto \Lambda_{Z^0} \oplus \Lambda_{T_G} \longto \pi_1( G) \longto 0
  \end{equation}
  is exact, because its subsequence $0 \longto \Lambda_{T_{\Gtilde}} \longto \Lambda_{T_G} \longto \pi_1( G) \longto 0$ is exact,
  and the resulting sequence of quotients $0 \longto \Lambda_{Z^0} = \Lambda_{Z^0} \longto 0 \longto 0$ is also exact.
  \begin{lemma}
    There is an exact sequence of reductive groups
    \begin{equation} \label{Ghat}
      1 \longto \Gtilde \longto \Ghat \longto[\dt] \Gm \longto 1
    \end{equation}
    and an extension $\pihat: \Ghat \longto G$ of $\pi: \Gtilde \longto G$ 
    such that $\pihat_*: \pi_1( \Ghat) \longto \pi_1( G)$ maps
    $1 \in \Z = \pi_1( \Gm) = \pi_1( \Ghat)$ to the given element $d \in \pi_1( G)$.
  \end{lemma}
  \begin{proof}
    We view the given $d \in \pi_1( G)$ as a coset $d \subseteq \Lambda_{T_G}$ modulo $\Lambda_{\coroots}$. Let
    \begin{equation*}
      \Lambda_{T_{\Ghat}} \subseteq \Lambda_{T_G} \oplus \Z
    \end{equation*}
    be generated by $\Lambda_{\coroots} \oplus 0$ and $(d, 1)$, and let
    \begin{equation*}
      (\pihat, \dt): \Ghat \longto G \times \Gm
    \end{equation*}
    be the reductive group with the same root system as $G$, whose maximal torus $T_{\Ghat} = \pihat^{-1}( T_G)$ has cocharacter lattice
    $\Hom( \Gm, T_{\Ghat}) = \Lambda_{T_{\Ghat}}$. As $\pi_*$ maps $\Lambda_{T_{\Gtilde}}$ isomorphically onto $\Lambda_{\coroots}$, we obtain an exact sequence
    \begin{equation*}
      0 \longto \Lambda_{T_{\Gtilde}} \longto[\pi_*] \Lambda_{T_{\Ghat}} \longto[\pr_2] \Z \longto 0,
    \end{equation*}
    which yields the required exact sequence \eqref{Ghat} of groups. By its construction, $\pihat_*$
    maps the canonical generator $1 \in \pi_1( \Gm) = \pi_1( \Ghat)$ to $d \in \pi_1( G)$.
  \end{proof}
  Let $\mu$ denote the kernel of the central isogeny $\zeta \cdot \pi: Z^0 \times \Gtilde \twoheadrightarrow G$. Then
  \begin{equation*}
    \psi: Z^0 \times \Ghat \longto G \times \Gm,
      \qquad (z^0, \ghat) \longmapsto \big( \zeta( z^0) \cdot \pihat( \ghat), \dt( \ghat) \big)
  \end{equation*}
  is by construction a central isogeny with kernel $\mu$.
  Hence the induced $1$--morphism
  \begin{equation*}
    \psi_*: \Bun^0_{Z^0} \times \Bun_{\Ghat}^1 \longto \Bun_G^d \times \Bun_{\Gm}^1
  \end{equation*}
  is faithfully flat by Lemma \ref{flat}. Restricting to the point
  $\Spec( k) \longto \Bun^1_{\Gm}$ given by a line bundle $L$ of
  degree $1$ on $C$, we get a faithfully flat $1$--morphism
  \begin{equation*}
    (\psi_*)_L: \Bun^0_{Z^0} \times \Bun_{\Ghat, L} \longto \Bun_G^d.
  \end{equation*}
  Since $\Gamma( \Bun^0_{Z^0} \times \Bun_{\Ghat, L}, \O) = k$ by 
  Proposition \ref{prop:twisted}(i) and Lemma \ref{lemma:Pic}(i),
  part (i) of the theorem follows. The group stack $\Bun_{\mu}$
  acts by tensor product on these two $1$--morphisms $\psi_*$
  and $(\psi_*)_L$, turning both into $\Bun_{\mu}$--torsors;
  cf. Example \ref{isogeny2torsor}. The idea is to descend
  line bundles along the torsor $(\psi_*)_L$.

  We choose a principal $T_{\Ghat}$--bundle $\xihat$ on $C$ together 
  with an isomorphism of line bundles $\dt_* \xihat \cong L$.
  Then $\xi := \pihat_*( \xihat)$ is a principal $T_G$--bundle on $C$; 
  their degrees $\deltahat := \deg( \xihat) \in \Lambda_{T_{\Ghat}}$
  and $\delta := \deg( \xi) \in \Lambda_{T_G}$ are lifts of $d \in \pi_1( G)$. The diagram
  \begin{equation} \label{group-square} \xymatrix{
    Z^0 \times T_{\Ghat} \ar[rr]^{\id \times \iota_{\Ghat}} \ar[d]^{\psi} && Z^0 \times \Ghat \ar[d]^{\psi}\\
    T_G \times \Gm \ar[rr]^{\iota_G \times \id} && G \times \Gm
  } \end{equation}
  of groups induces the right square in the $2$--commutative diagram
  \begin{equation} \label{pihat_times_dt-square} \xymatrix{
    \Bun^0_{Z^0} \times \Bun^0_{T_{\Ghat}} \ar[rr]^{\id \times t_{\xihat}} \ar[d]^{\psi_*}
      && \Bun^0_{Z^0} \times \Bun_{T_{\Ghat}}^{\deltahat} \ar[rr]^{(\id \times \iota_{\Ghat})_*} \ar[d]^{\psi_*} && \Bun^0_{Z^0} \times \Bun_{\Ghat}^1 \ar[d]^{\psi_*}\\
    \Bun^0_{T_G} \times \Bun_{\Gm}^1 \ar[rr]^{t_{\xi} \times \id} && \Bun_{T_G}^{\delta} \times \Bun_{\Gm}^1 \ar[rr]^{(\iota_G \times \id)_*}
      && \Bun_G^d \times \Bun_{\Gm}^1
  } \end{equation}
  of moduli stacks; note that $t_{\xihat}$ and $t_{\xi}$ are equivalences.
  Restricting the outer rectangle again to the point $\Spec( k) \longto \Bun^1_{\Gm}$
  given by $L$, we get the diagram
  \begin{equation} \label{pihat-square} \xymatrix{
    \Bun^0_{Z^0 \times T_{\Gtilde}} \ar@{=}[r]^-{\sim} \ar[dr]_{(\zeta \cdot \pi)_*}
      & \Bun^0_{Z^0} \times \Bun^0_{T_{\Gtilde}} \ar[rr]^{\id \times \iota_{\xihat}} \ar[d]^{\zeta_* \otimes \pi_*} && \Bun^0_{Z^0} \times \Bun_{\Ghat, L} \ar[d]^{(\psi_*)_L}\\
    & \Bun^0_{T_G} \ar[rr]^{(\iota_G)_* \circ t_{\xi}} && \Bun_G^d
  } \end{equation}
  containing an instance $\iota_{\xihat}$ of the $1$--morphism 
  \eqref{iota_xi} defined in Subsection \ref{sc2torus}.
  According to the Proposition \ref{prop:torus} and 
  Proposition \ref{prop:twisted},
  \begin{equation*}
    \iota_{\xihat}^*: \relPic( \Bun_{\Ghat, L}) \longto \relPic( \Bun^0_{T_{\Gtilde}})
  \end{equation*}
  is a morphism of group schemes over $k$. This morphism is a closed 
  immersion, according to Proposition \ref{prop:twisted}(iii),
  if $g_C \geq 1$ or if $\xihat$ is chosen appropriately, as
  explained in Lemma \ref{injective}; we assume this in the sequel.
  Using Lemma \ref{see-saw} and Corollary \ref{torus:product}, it 
  follows that
  \begin{equation*}
    (\id \times \iota_{\xihat})^*:
      \relPic( \Bun^0_{Z^0}) \oplus \relPic( \Bun_{\Ghat, L}) \cong \relPic( \Bun^0_{Z^0} \times \Bun_{\Ghat, L}) \longto \relPic( \Bun^0_{Z^0 \times T_{\Gtilde}})
  \end{equation*}
  is a closed immersion of group schemes over $k$ as well.

  The group stack $\Bun_{\mu}$ still acts by tensor product on the 
  vertical $1$--morphisms in \eqref{pihat_times_dt-square} and in \eqref{pihat-square}.
  Since the diagram \eqref{group-square} of groups is cartesian, 
  \eqref{pihat_times_dt-square} and \eqref{pihat-square} are morphisms of 
  $\calG$--torsors; cf. Example \ref{pullback2torsors}. Proposition 
  \ref{torsors:cartesian} applies to the latter morphism of torsors,
  yielding a cartesian square
  \begin{equation} \label{Pics:cartesian} \xymatrix{
    \relPic( \Bun_G^d) \ar[d]^{\psi^*_L} \ar[rr]^-{t_{\xi}^* \circ \iota_G^*} && \relPic( \Bun^0_{T_G}) \ar[d]^{(\zeta \cdot \pi)^*}\\
    \relPic( \Bun^0_{Z^0}) \oplus \relPic( \Bun_{\Ghat, L}) \ar[rr]^-{(\id \times \iota_{\xihat})^*} && \relPic( \Bun^0_{Z^0 \times T_{\Gtilde}})
  } \end{equation}
  of Picard functors. Thus $\relPic( \Bun_G^d)$ is representable, and 
  $t_{\xi}^* \circ \iota_G^*$ is a closed immersion;
  this proves part (ii) of the theorem.

  The image of the mapping cone \eqref{lattices:exact} under the
  exact functor $\relHom( \_, J_C)$, and the mapping cones of the
  two cartesian squares given by diagram \eqref{Pics:cartesian}
  and Lemma \ref{NS:cartesian}, are the columns of the commutative diagram
  \begin{equation*} \xymatrix{
    & 0 \ar[d] & 0 \ar[d] && 0 \ar[d] & \\
    & \relHom( \pi_1( G), J_C) \ar[d] & \relPic( \Bun_G^d) \ar[d] && \NS( \Bun_G^d) \ar[d] & \\
    0 \ar[r] & *\txt{$\relHom( \Lambda_{Z^0}, J_C)$\\$\oplus$\\$\relHom( \Lambda_{T_G}, J_C)$} \ar[r]^-{j_{Z^0} \oplus j_{T_G}} \ar[d]
      & *\txt{$\relPic( \Bun_{Z^0}^0)$\\$\oplus$\\$\relPic( \Bun_{\Ghat, L})$\\$\oplus$\\$\relPic( \Bun_{T_G}^0)$} \ar[d]
        \ar[rr]^-{c_{Z^0} \oplus c_{\Gtilde} \oplus c_{T_G}}
      && *\txt{$\NS(     \Bun_{Z^0})  $\\$\oplus$\\$\NS(     \Bun_{\Gtilde}) $\\$\oplus$\\$\NS(     \Bun_{T_G})$} \ar[d] \ar[r] & 0\\
    0 \ar[r] & \relHom( \Lambda_{Z^0 \times T_{\Gtilde}}, J_C) \ar[d] \ar[r]^-{j_{Z^0 \times T_{\Gtilde}}}
      & \relPic( \Bun^0_{Z^0 \times T_{\Gtilde}}) \ar[rr]^-{c_{Z^0 \times T_{\Gtilde}}} && \NS( \Bun_{Z^0 \times T_{\Gtilde}}) \ar[r] & 0\\
    & 0 &&&&
  } \end{equation*}
  whose two rows are exact due to Proposition \ref{prop:torus}(ii) and 
  Proposition \ref{prop:twisted}(ii). Applying the snake lemma to
  this diagram, we get an exact sequence
  \begin{equation*}
    0 \longto \relHom( \pi_1( G), J_C) \xrightarrow{ j_G( \iota_G, \delta)} \relPic( \Bun_G^d) \xrightarrow{ c_G( \iota_G, \delta)} \NS( \Bun_G^d) \longto 0.
  \end{equation*}
  The image of $j_G( \iota_G, \delta)$ and the kernel of $c_G( \iota_G, \delta)$
  are a priori independent of the choices made, since both are the
  largest quasicompact open subgroup in $\relPic( \Bun_G^d)$. If $G$ is 
  a torus and $d = 0$, then this is the exact sequence of Proposition
  \ref{prop:torus}; in general, the construction provides a morphism of exact sequences
  \begin{equation} \label{embed2torus} \xymatrix{
    0 \ar[r] & \relHom( \pi_1( G), J_C) \ar[rr]^-{j_G( \iota_G, \delta)} \ar[d]^{\pr^*}
      && \relPic( \Bun_G^d) \ar[d]^{t_{\xi}^* \circ \iota_G^*} \ar[rr]^{c_G( \iota_G, \delta)} && \NS( \Bun_G^d) \ar[r] \ar[d]^{(\iota_G)^{\NS, \delta}} & 0\\
    0 \ar[r] & \relHom( \Lambda_{T_G}, J_C) \ar[rr]^-{j_{T_G}} && \relPic( \Bun^0_{T_G}) \ar[rr]^{c_{T_G}} && \NS( \Bun_{T_G}) \ar[r] & 0
  } \end{equation}
  whose three vertical maps are all injective. Using Proposition 
  \ref{prop:torus}(iii), this implies that $j_G( \iota_G, \delta)$ and 
  $c_G( \iota_G, \delta)$ depend at most on the choice of
  $\iota_G: T_G \hookrightarrow G$ and of $\delta$, but not on
  the choice of $\Ghat$, $L$ or $\xihat$; thus the notation.
  Together with the following two lemmas, this proves
  the remaining parts (iii) and (iv) of the theorem.
\end{proof}
\begin{lemma}
  The above map $j_G( \iota_G, \delta): \relHom( \pi_1( G), J_C) \longto \relPic( \Bun_G^d)$
  \begin{itemize}
   \item[i)] does not depend on the lift $\delta \in \Lambda_{T_G}$ of $d \in \pi_1( G)$,
   \item[ii)] does not depend on the maximal torus $\iota_G: T_G \hookrightarrow G$, and
   \item[iii)] satisfies $\varphi^* \circ j_H = j_G \circ \varphi^*: \relHom( \pi_1( H), J_C) \longto \relPic( \Bun_G^d)$
    for all $\varphi: G \longto H$.
  \end{itemize}
\end{lemma}
\begin{proof}
  If $G$ is a torus, then $\delta$ and $\iota_G$ are unique,
  so (i) and (ii) hold trivially.

  The claim is empty for $g_C = 0$, so we assume $g_C \geq 1$. Then
  the above construction works for all lifts $\delta$ of $d$, because
  $\iota_{\xihat}^*$ is a closed immersion for all $\xihat$.

  Given $\varphi: G \longto H$ and a maximal torus
  $\iota_H: T_H \hookrightarrow H$ with $\varphi( T_G) \subseteq T_H$,
  we again put $e := \varphi_* d \in \pi_1( H)$ and
  $\eta := \varphi_* \delta \in \Lambda_{T_H}$. Then the diagram
  \begin{equation} \label{j_G:functorial} \xymatrix{
    \relHom( \pi_1( H), J_C) \ar[rr]^-{j_H( \iota_H, \eta)} \ar[d]^{\varphi^*} && \relPic( \Bun_H^e) \ar[d]^{\varphi^*}\\
    \relHom( \pi_1( G), J_C) \ar[rr]^-{j_G( \iota_G, \delta)}                 && \relPic( \Bun_G^d)
  } \end{equation}
  commutes, because it commutes after composition with the closed immersion
  \begin{equation*}
    t_{\xi}^* \circ \iota_G^*: \relPic( \Bun_G^d) \longto \relPic( \Bun^0_{T_G})
  \end{equation*}
  from diagram \eqref{embed2torus}, using Remark \ref{c_T:functorial}. 
  In particular, (iii) follows from (i) and (ii).

  i) For $G = \GL_2$, it suffices to take $\varphi = \det: \GL_2 \longto \Gm$
  in the above diagram \eqref{j_G:functorial}, since
  $\det_*: \pi_1( \GL_2) \longto \pi_1( \Gm)$ is an isomorphism.

  For $G = \PGL_2$, it then suffices to take $\varphi = \pr: \GL_2 \twoheadrightarrow \PGL_2$ in the same diagram \eqref{j_G:functorial},
  since $\pr_*: \pi_1( \GL_2) \longto \pi_1( \PGL_2)$ is surjective.

  As (i) holds trivially for $G = \SL_2$, and clearly holds for $G \times \Gm$
  if it holds for $G$, this proves (i) for all groups $G$ of semisimple rank one.

  In the general case, let $\alpha^{\vee} \in \Lambda_{T_G}$ be a coroot,
  and let $\varphi: G_{\alpha} \hookrightarrow G$ be the corresponding subgroup of
  semisimple rank one. Then the diagram \eqref{j_G:functorial} shows
  $j_G( \iota_G, \delta) = j_G( \iota_G, \delta + \alpha^{\vee})$, since
  $\varphi_*: \pi_1( G_{\alpha}) \longto \pi_1( G)$ is surjective. This 
  completes the proof of i, because any two lifts $\delta$ of $d$
  differ by a sum of coroots.

  ii) now follows from Weyl--invariance; cf. Subsection \ref{sc:NS}.
\end{proof}
\begin{lemma}
  The above map $c_G( \iota_G, \delta): \relPic( \Bun_G^d) \longto \NS( \Bun_G^d)$
  \begin{itemize}
   \item[i)] does not depend on the lift $\delta \in \Lambda_{T_G}$ of $d \in \pi_1( G)$,
   \item[ii)] does not depend on the maximal torus $\iota_G: T_G \hookrightarrow G$, and
   \item[iii)] satisfies $\varphi^{\NS, d} \circ c_H = c_G \circ \varphi^*: \relPic( \Bun_H^e) \longto \NS( \Bun_G^d)$
    for all $\varphi: G \longto H$.
  \end{itemize}
\end{lemma}
\begin{proof}
  If $G$ is a torus, then $\delta$ and $\iota_G$ are unique; if $G$ is simply
  connected, then $c_G( \iota_G, \delta)$ coincides by construction with the
  isomorphism $c_G$ of Proposition \ref{prop:twisted}(ii). In both cases,
  (i) and (ii) follow, and we can use the notation $c_G$ without ambiguity. 

  Given a representation $\rho: G \longto \SL( V)$, the diagram
  \begin{equation} \label{c_G:representation} \xymatrix{
    \relPic( \Bun_{\SL( V)}) \ar[rr]^-{c_{\SL( V)}} \ar[d]^{\rho^*} && \NS( \Bun_{\SL( V)}) \ar[d]^{\rho^{\NS, d}}\\
    \relPic( \Bun_G^d) \ar[rr]^-{c_G( \iota_G, \delta)}             && \NS( \Bun_G^d)
  } \end{equation}
  commutes, because it commutes after composition with the injective map
  \begin{equation*}
    (\iota_G)^{\NS, \delta}: \NS( \Bun_G^d) \longto \NS( \Bun_{T_G})
  \end{equation*}
  from diagram \eqref{embed2torus}, using Lemma 
  \ref{iota^NS_functorial}, Corollary \ref{cor:sl}, Remark 
  \ref{c_T:functorial}, and the $2$--commutative squares
  \begin{equation*} \xymatrix{
    \Bun^0_{T_G} \ar[r]^{t_{\xi}} \ar[d]^{\rho_*} & \Bun^{\delta}_{T_G} \ar[r]^{(\iota_G)_*} \ar[d]^{\rho_*} & \Bun_G^d \ar[d]^{\rho_*}\\
    \Bun^0_{T_{\SL( V)}} \ar[r]^{t_{\rho_* \xi}} & \Bun^{\rho_* \delta}_{T_{\SL( V)}} \ar[r]^{\iota_*} & \Bun_{\SL( V)})
  } \end{equation*} 
  in which $\iota: T_{\SL( V)} \hookrightarrow \SL( V)$ is a maximal torus containing $\rho( T_G)$.

  Similarly, given a homomorphism $\chi: G \longto T$ to a torus $T$, 
  the diagram
  \begin{equation} \label{c_G:character} \xymatrix{
    \relPic( \Bun_T^{\chi_* d}) \ar[rr]^-{c_T} \ar[d]^{\chi^*} && \NS( \Bun_T) \ar[d]^{\chi^*}\\
    \relPic( \Bun_G^d) \ar[rr]^-{c_G( \iota_G, \delta)}             && \NS( \Bun_G^d)
  } \end{equation}
  commutes, again because it commutes after composition with the same
  injective map $(\iota_G)^{\NS, \delta}$ from diagram \eqref{embed2torus},
  using Lemma \ref{iota^NS_functorial}, Remark \ref{c_T:functorial}, and 
  the $2$--commutative squares
  \begin{equation*} \xymatrix{
    \Bun^0_{T_G} \ar[r]^{t_{\xi}} \ar[d]^{\chi_*} & \Bun^{\delta}_{T_G} \ar[r]^{(\iota_G)_*} \ar[d]^{\chi_*} & \Bun_G^d \ar[d]^{\chi_*}\\
    \Bun^0_T \ar[r]^{t_{\chi_* \xi}} & \Bun^{\chi_* \delta}_T \ar@{=}[r] & \Bun^{\chi_* \delta}_T.
  } \end{equation*} 

  The two commutative diagrams \eqref{c_G:representation} and \eqref{c_G:character}
  show that the restriction of $c_G( \iota_G, \delta)$ to the images of all $\rho^*$
  and all $\chi^*$ in $\relPic( \Bun_G^d)$ modulo $\relHom( \pi_1( G), J_C)$ does
  not depend on the choice of $\delta$ or $\iota_G$. But these images generate a
  subgroup of finite index, according to Proposition \ref{NS_extension} and Remark
  \ref{simply_connected_products}. Thus (i) and (ii) follow. The functoriality in
  (iii) is proved similarly; it suffices to apply these arguments to homomorphisms
  $\rho: H \longto \SL( V)$, $\chi: H \longto T$ and their compositions with
  $\varphi: G \longto H$, using Corollary \ref{phi^NS_functorial}.
\end{proof}

\end{document}